\documentstyle[amssymb,12pt]{amsart}
\input xy
\xyoption{all}

\newdir{((}{{\lhook\kern-.1em}}                                              
\newdir{))}{{\rhook\kern-.1em}}                                              
\newdir{ >}{{}*!/-5pt/@{>}}                                                  

\makeatletter
\newcommand{\subsectionr}
{\@startsection{subsection}{3}{0pt}{\baselineskip}
 {-\fontdimen2\font}{\normalfont\bfseries}}
\newcommand{\subsubsectionr}
{\@startsection{subsubsection}{3}{0pt}{\baselineskip}
 {-\fontdimen2\font}{\normalfont\bfseries}}
\renewcommand{\subsubsection}
{\@startsection{subsubsection}{3}{0pt}{\baselineskip}
 {0.05\baselineskip}{\normalfont\bfseries}}
\makeatother

\newcommand{\la}{\langle}
\newcommand{\ra}{\rangle}

\newtheorem{theorem}{\bf Theorem}[section]
\newtheorem{lemma}[theorem]{\bf Lemma}

\newtheorem{corollary}[theorem]{\bf Corollary}

{%\theorembodyfont{\rmfamily}
\newtheorem{definition}[theorem]{\bf Definition}

}

\newcommand{\be}{\begin{equation}}             
\newcommand{\ee}{\end{equation}}        

\newfont{\bfc}{cmbsy10 scaled 1200}  % bold face calligraphic
\newfont{\dr}{msbm10 scaled \magstep1}  %letra doble raya
\newfont{\sdr}{msbm8}  % small letra doble raya
\newfont{\gl}{eufm10 scaled \magstep1}  % german letters as well

\DeclareFontFamily{OT1}{rsfs}{}
\DeclareFontShape{OT1}{rsfs}{n}{it}{<->rsfs10}{}
\DeclareMathAlphabet{\curly}{OT1}{rsfs}{n}{it}

\newcommand{\CC}{{\Bbb C}}
\newcommand{\CP}{{\Bbb CP}}

\newcommand{\EE}{{\Bbb E}}

\newcommand{\KK}{{\mathbf K}}

\newcommand{\NN}{{\Bbb N}}

\newcommand{\RR}{{\Bbb R}}
\renewcommand{\SS}{{\mathbf S}}
\newcommand{\TT}{{\mathbf T}}
\newcommand{\UU}{{\mathbf U}}

\newcommand{\ZZ}{{\Bbb Z}}

\newcommand{\ov}{\overline}

\newcommand{\fX}{{\curly X}}

\newcommand{\Ad}{\operatorname{Ad}}
\newcommand{\Aut}{\operatorname{Aut}}

\newcommand{\codim}{\operatorname{codim}}

\newcommand{\Coker}{\operatorname{Coker}}
\newcommand{\cusp}{\operatorname{c}}

\newcommand{\diam}{\operatorname{diam}}
\newcommand{\dist}{\operatorname{dist}}

\newcommand{\End}{\operatorname{End}}
\newcommand{\ev}{\operatorname{ev}}
\newcommand{\fibr}{\operatorname{fibr}}
\newcommand{\Fl}{\operatorname{Fl}}
\newcommand{\GL}{\operatorname{GL}}

\newcommand{\Hom}{\operatorname{Hom}}
\newcommand{\Id}{\operatorname{Id}}

\newcommand{\inter}{\operatorname{int}}
\newcommand{\Ind}{\operatorname{Ind}}
\newcommand{\Jac}{\operatorname{Jac}}
\newcommand{\Ker}{\operatorname{Ker}}
\newcommand{\Lie}{\operatorname{Lie}}

\newcommand{\Map}{\operatorname{Map}}
\newcommand{\op}{\operatorname{op}}

\newcommand{\reg}{\operatorname{reg}}
\newcommand{\rk}{\operatorname{rk}}
\newcommand{\PSL}{\operatorname{PSL}}

\newcommand{\Tr}{\operatorname{Tr}}

\newcommand{\Vol}{\operatorname{Vol}}
\renewcommand{\exp}{\operatorname{exp}}

\newcommand{\AAA}{{\curly A}}
\newcommand{\BBB}{{\curly B}}
\newcommand{\CCC}{{\curly C}}
\newcommand{\DDD}{{\cal D}}
\newcommand{\EEE}{{\cal E}}
\newcommand{\FFF}{{\cal F}}
\newcommand{\GGG}{{\curly G}}
\newcommand{\HHH}{{\cal H}}
\newcommand{\III}{{\curly I}}

\newcommand{\MMM}{{\cal M}}

\newcommand{\PPP}{{\curly P}}
\newcommand{\QQQ}{{\curly Q}}
\newcommand{\RRR}{{\cal R}}
\newcommand{\SSS}{{\curly S}}

\newcommand{\WWW}{{\cal W}}
\newcommand{\YMH}{{\cal YMH}}

\newcommand{\imag}{{\mathbf i}}
\newcommand{\qu}{/\kern-.7ex/}
\newcommand{\exh}{\to\kern-1.8ex\to}
\newcommand{\nsubset}{\subset\kern-2ex/\kern+1ex}
\newcommand{\bM}{\widetilde{\MMM}}
\newcommand{\bW}{\widetilde{\WWW}}
\newcommand{\bR}{\widetilde{\RRR}}

\newcommand{\bE}{\widetilde{E}}
\newcommand{\bEE}{\widetilde{\EE}}

\renewcommand{\Im}{\operatorname{Im}}

\newcommand{\bPhi}{{\mathbf \Phi}}                 

\newcommand{\obPhi}{\widetilde{\bPhi}}

\title{Hamiltonian Gromov--Witten invariants}
\author{Ignasi Mundet i Riera}
\date{28--1--2000}
\address{Centre de Math{\'e}matiques CNRS UMR 7640, 
{\'E}cole Polytechnique, Palaiseau,
France, and Departamento de Matem{\'a}ticas,
Universidad Aut{\'o}noma de Madrid,
Madrid, Spain}
\email{ignasi@@math.polytechnique.fr}
 % MSC2000

\begin{document}
\maketitle

\begin{abstract}
In this paper we introduce invariants of semi-free
Hamiltonian actions of $S^1$ on compact symplectic manifolds (which
satisfy some technical conditions related to positivity)
using the space of solutions to certain gauge theoretical
equations. These equations generalize at the same time
the vortex equations and the holomorphicity equation used in 
Gromov--Witten theory. In the definition of the invariants we combine
ideas coming from gauge theory and the ideas underlying the
construction of Gromov--Witten invariants. 
\end{abstract}

\tableofcontents

\section{Introduction}
\label{s:introduction}
Gromov--Witten invariants are among the most important and useful
tools in symplectic topology. They can be used for example to prove 
the existence of two symplectic structures on a compact smooth
manifold which are not deformation equivalent (see \cite{Ru1}).
Other applications appear in the study of 
symplectic fibrations (see \cite{McD2}), the topology of the group
of Hamiltonian symplectomorphisms (see \cite{Se}), and Weinstein's
conjecture \cite{LiuTi}.
Finally, Gromov--Witten invariants play a prominent role in the
celebrated mirror conjecture (see for example \cite{Gi}).

It is a natural question whether Gromov--Witten invariants or some
related construction can be
used to study Hamiltonian actions of compact Lie groups
on symplectic manifolds. The purpose of this paper is to make some
steps towards an
affirmative answer to this question. 
For any symplectic manifold with a Hamiltonian action of $S^1$ 
we introduce a set of equations whose moduli space of solutions, 
in a way very much similar to Gromov--Witten theory, allow to define
invariants of the symplectic manifold together with the Hamiltonian
action. 
These equations generalise on the one hand the holomorphicity
equation used in Gromov--Witten theory and on the other hand the
gauge theoretical
vortex equations or their analogues as considered in full
generality by Banfield \cite{Ba}. Although we restrict
ourselves to Hamiltonian actions of $S^1$, the definition of the
invariants can be given (at least heuristically) for any action of a
compact group.
We call these new invariants Hamiltonian Gromov--Witten invariants.

Let $(F^{2n},\omega)$ be a compact symplectic manifold supporting a 
Hamiltonian
action of $S^1$ with moment map $\mu:F\to\Lie(S^1)^*=(\imag\RR)^*$.
Let $ES^1\to BS^1$ be the universal principal $S^1$-bundle.
Recall that the Borel construction of $F$ is 
$$F_{S^1}=ES^1\times_{S^1}F,$$
and that the equivariant homology (resp. cohomology) of $F$ is by
definition $H_*^{S^1}(F;\ZZ):=H_*(F_{S^1};\ZZ)$ 
(resp. $H^*_{S^1}(F;\ZZ):=H^*(F_{S^1};\ZZ)$).
The Hamiltonian Gromov--Witten invariants depend on the choice of an
element in $H_2^{S^1}(F;\ZZ)$, a collection of elements in
$H^*_{S^1}(F;\ZZ)$, and an element of $\Lie(S^1)$. 
It is an interesting question to relate these invariants to the
equivariant Gromov--Witten invariants \cite{GiKm,Lu}
(note, however, that in the definition of
equivariant Gromov--Witten invariants one does not use any element of
$\Lie(S^1)$, and that the action of $S^1$ does not need to be
Hamiltonian for these invariants to be defined).

To define the Hamiltonian Gromov--Witten invariants we follow the
ideas and techniques used in the definition of Gromov--Witten
invariants in \cite{McDS1}. We are forced to assume some technical
conditions on the manifold $F$ and on the action. However, we expect
that using the techniques of virtual moduli cycles (as developped in
the several papers which give a construction of Gromov--Witten
invariants for general compact symplectic manifolds, see 
\cite{FuOn,LiTi,Ru2,Si}) 
a definition of Hamiltonian Gromov--Witten invariants could be given
in full generality. 

This paper is essentially based on the Ph.D. Thesis of the author
\cite{Mu2}, which was submitted in the Universidad Aut{\'o}noma de Madrid
in the spring of 1999. 
The author discovered the equations which are used to define the
Hamiltonian Gromov--Witten invariants inspired by previous work on
Hitchin--Kobayashi correspondence (see \cite{Mu1}).
In the summer of 1999 the author knew that
K. Cieliebak, A.R. Gaio and D. Salamon
had independently discovered the same equations and how to define the 
invariants (see \cite{CiGaSa,Ga}). On the other hand, A. Bertram,
G. Daskalopoulos and R. Wentworth studied in \cite{BtDaWe}  
Gromov--Witten invariants of Grassmannians using ideas similar to
the ones which we use.

This paper is organised as follows. 
In the rest of this section we introduce the equations and we give
a heuristic definition of the invariants. In Section
\ref{s:smoothness} we introduce Sobolev completions of our parameter
space and we explain how to perturb the equations in order to
get a smooth moduli space of solutions. We also prove a regularity
result for Sobolev solutions of the equations.
In Section \ref{s:compactification} we define a compactification of
the moduli space of solutions to the equation. This compactification
generalises Gromov compactification of the moduli space of
pseudo-holomorphic curves. In Section \ref{s:rationalcurves} we study
the moduli of rational curves for a generic $S^1$-invariant complex
structure on $F$. Finally, in Section \ref{s:definition} we recall the
basic definition of the theory of pseudo-cycles and we use them to
give a rigorous definition of the invariants under the conditions
specified in Subsection \ref{conditions}.

\noindent{\bf Acknowledgements.}
I am very much indebted to my thesis advisor, Oscar Garc{\'\i}a--Prada, for his
continuous support, encouragement and generous 
share of ideas. I also would like to thank G. Segal for a conversation
which helped me to understand the topological construction underlying
the definition of the invariants.

\subsection{The equations}
\label{theequations}
Take a $S^1$-invariant complex structure $I\in\End(TF)$ such that 
$g(\cdot,\cdot)=\omega(\cdot,I\cdot)$ is a Riemannian metric on $F$.
Such complex structures always exists, and they form a contractible space
(see Lemma 5.49 in \cite{McDS2}); this implies that the Chern classes
of the complex bundle $(TF,I)$ only depend on (the deformation class
of) $\omega$.

Let $\Sigma$ be a compact connected Riemann surface, 
with a fixed Riemannian metric.
Let $E\to\Sigma$ be a principal $S^1$-bundle, and let 
$\pi:\FFF=E\times_{S^1}F\to\Sigma$ be the associated fibration with fibre $F$.
Since the moment map is by definition $S^1$-invariant, it extends to give
a map $\mu:\FFF\to(\imag\RR)^*$. 
Let $\AAA$ be the space of connections on $E$, 
let $\GGG=\Map(\Sigma,S^1)$ be the gauge group of $E$,
and let $\SSS=\Gamma(\FFF)$ be the space of smooth sections of $\FFF$.

Any connection $A\in\AAA$ induces a projection 
$\alpha_A:T\FFF\to T\FFF_v=\Ker d\pi$. Using this map we define the
covariant derivative with respect to $A$ of a section $\phi\in\SSS$ to be
$$d_A\phi=\alpha_A\circ d\phi\in\Omega^1(\Sigma;\phi^* T\FFF_v).$$
On the other hand, since $I$ is $S^1$-invariant it can be extended to
give a complex structure on $T\FFF_v$. Hence, we can split $d_A\phi$ as
the sum of its holomorphic part 
$\partial_A\phi\in\Omega^{1,0}(\Sigma;\phi^* T\FFF_v)$
plus its antiholomorphic part
$\ov{\partial}_A\phi\in\Omega^{0,1}(\Sigma;\phi^* T\FFF_v)$.

Let $\Lambda:\Omega^2(\Sigma)\to\Omega^0(\Sigma)$ denote the contraction
with the volume form on $\Sigma$. We will henceforth identify 
$\Lie(S^1)^*=(\imag\RR)^*$ with $\Lie(S^1)=(\imag\RR)$
by assigning to $a\in\imag\RR$ the element $a^*:\imag\RR\to\RR$ which maps
$b\in\imag\RR$ to $a^*(b)=\la a,b\ra=-ab\in\RR$.
Let $c\in\imag\RR$.
We will consider the following equations on $(A,\phi)\in\AAA\times\SSS$
\begin{equation}
\left\{\begin{array}{l}
\ov{\partial}_A\phi=0\\
\Lambda F_A+\mu(\phi)=c
\end{array}\right.
\label{equs}
\end{equation}
where $F_A\in\Omega^2(\Sigma;\imag\RR)$ denotes the curvature of $A$.
Taking on $\AAA\times\SSS$ the diagonal action of the gauge group (acting
by pullback both in $\AAA$ and $\SSS$), the set of solutions
to the equations is $\GGG$-invariant.

\begin{definition}
We call a tuple $(\Sigma,E,A,\phi,c)$ for which (\ref{equs}) is
satisfied a twisted holomorphic curve, THC for short.
\end{definition}

\noindent{\bf Example.}
When $F=\CC$ and $S^1\subset\CC$ acts on it by multiplication,
the equations (\ref{equs}) coincide with the abelian vortex equations (see
for example \cite{GP}). (Note, however, that in this paper we assume
that $F$ is compact.)

\subsection{Relation with holomorphic curves in $\FFF$}
\label{corbeshol}
Any connection $A\in\AAA$ gives rise to a splitting 
$T\FFF\simeq T\FFF_v\oplus \pi^*T\Sigma$ (which is given by
taking as complementary to $T\FFF_v\subset T\FFF$ the kernel 
$\Ker\alpha_A$). Using this splitting we may combine the complex
structure on $T\FFF_v$ (which was induced by the complex structure on
$F$) with the complex structure of $T\Sigma$ to get a complex structre
$I(A)\in\End(T\FFF)$. It is then straightforward to prove the
following.

\begin{lemma}
\label{lcorbeshol}
A section $\phi\in\SSS$ is $I(A)$ holomorphic as a map from $\Sigma$
to $\FFF$ if and only if $\ov{\partial}_A\phi=0$.
\end{lemma}

\subsection{Principal bundles and maps to Borel construction}
\label{toprest}
In this subsection we will describe a construction which relates
principal bundles on a topological space and sections of an associated
bundle to maps into the Borel construction of the fibre of the
associated bundle. We will use this construction twice in this paper:
first, to fix the topology of the bundle $E\to\Sigma$ and of the section
$\phi\in\Gamma(E\times_{S^1}F)$; and, later, to build an {\it
  equivariant evaluation map} which will play the role of the
evaluation map in Gromov--Witten theory.

Let $M$ be a CW-complex, let $K$ be a compact connected Lie group, 
and let $V$ be a CW-complex with a left action of $K$. 
Let $EK\to BK$ be the universal principal $K$-bundle\footnote{We 
assume that the action of the structure group 
on any principal bundle is on the right.}. We denote by
$V_K:=EK\times_K V$ the Borel construction of $V$.

\begin{lemma}
Let $\GGG_P=\Gamma(P\times_{\Ad}K)$ be the gauge group of $P$.
There is a canonical bijection between the set of homotopy classes of
maps $[M,V_K]$ and the set of homotopy classes of pairs
consisting of a principal $K$-bundle $P\to M$  
and a gauge equivalence class $\GGG_P\sigma\subset\Gamma(P\times_KV)$ 
of sections of the associated bundle $P\times_KV$.
\label{bijeccio}
\end{lemma}
\begin{pf}
Take a pair $(P,\sigma)$, where $P\to M$ is a principal $K$-bundle and
$\sigma\in\Gamma(P\times_KV)$. 
Let $C$ be the set of $K$-equivariant continuous maps $P\to
EK$. Note that any element of $C$ gives a lift of the classifying map
$M\to BK$ of $P$. Denote by $EK^{\op}$ the space $EK$ with the action
$\rho_L:K\times EK\to EK$ of $K$ on the left given by 
$$\rho_L(k,x)=\rho_R(x,k^{-1}),$$
where $\rho_R:EK\times K\to EK$ is the usual right action. There is a
canonical bijection $$C\simeq\Gamma(P\times_K EK^{\op}),$$
where $\Gamma$ denotes the space of continuous sections. Since
$EK^{\op}$ is contractible, we deduce that $C$ is nonempty and
contractible. Fix one such map $C\ni c:P\to EK$. 
If $g\in\GGG_P$, then $c\circ g$ is also $K$-equivariant, and
it is homotopic (as $K$-equivariant maps) to $c$:
\begin{equation}
c\circ g\sim c.
\label{eshomotopic}
\end{equation}
The section $\sigma$ gives a $K$-equivariant map $\psi:P\to V$ 
(i.e. $\psi(pk)=k^{-1}\psi(p)$). Hence the map
$(c,\psi):P\to EK\times V$ descends to give
a map $\sigma_P:M\to V_K$. 
If $g\in\GGG_K$ then the section $g^*\sigma$ induces the $K$-equivariant
map $\psi\circ g:P\to V$. Now, by (\ref{eshomotopic})
$$(c,\psi\circ g)\sim (c\circ g^{-1},\psi\circ g),$$
and hence $(g^*\sigma)_P\simeq\sigma_P$. Finally, it is clear that the
homotopy class $[\sigma_P]$ only depends on the homotopy class of the
section $\sigma$.

Conversely, if $f: M\to V_K$ is any continuous
map, we let $P_f$ be the pullback $f^*\pi_K^*EK$, where
$\pi_K:V_K\to BK$ is the projection. 
Since there is a canonical isomorphism
$P_f\times_KV\simeq f^*\pi_K^*V_K$, the map $f$
induces a section $\sigma_f\in\Gamma(P_f\times_KV)$. This 
construction is the inverse of the preceeding one. In particular, we
have
\begin{equation}
P\simeq \sigma_f^*\pi_K^* EK.
\label{recupera}
\end{equation}
\end{pf}

\noindent{\bf Remark.} Another way to look at this lemma is the
following. Given a $K$ principal bundle $P\to M$ there is, up to
homotopy, a unique $K$-equivariant map $P\to EK$. Combining it with
the identity map from $V$ to $V$ and quotienting by $K$ we get a
map $P\times_KV\to V_K$ which induces a canonical map in cohomology
$c_P^*:H^*_K(V;\ZZ)\to H^*(P\times_K V;\ZZ)$. 
Furthermore, the image of $c_P^*$ lies inside the fixed point set of the action
of $\GGG_P$ on $H^*(P\times_KV;\ZZ)$.
Then, if $\sigma\in\Gamma(P\times_KV)$, we have 
$$\sigma_P^*=\sigma^*c_P^*.$$
An important consequence of this equality is the following. If 
$$\xymatrix{P\ar[r]^{s}\ar[d] & P' \ar[d]\\ M\ar[r] & M'}$$
is a Cartesian diagram of $K$ principal bundles, 
and if we denote $s_V:P\times_KV\to
P'\times_KV$ the map induced by $s$, then for any section
$\sigma\in\Gamma(P\times_KV)$ we have
\begin{equation}
\sigma_P^*=(s_V\sigma)^* c_{P'}.
\label{holabondia}
\end{equation}

\subsection{The moduli space}
Similarly to what one does in Gromov--Witten theory, we will study the set of
solutions to equations (\ref{equs}) for all pairs $(E,\phi)$
satisfying certain topological constraints. In Gromov--Witten theory
the constraints are specified by the choice of an element of $H_2$ of
the target space. In our case we will need to chose an element of
$H_2^{S_1}(F;\ZZ)$. 
We will apply the construction of Lemma \ref{bijeccio} to $S^1$
principal bundles $E\to\Sigma$. Observe that 
if $\sigma\in\Gamma(E\times_{S^1}F)$, 
we may recover the degree of $E$ using $\sigma_E$, since by
(\ref{recupera}) 
$$\deg E=\la c_1(E),[\Sigma]\ra
=\la c_1(ES^1),{\pi_{S^1}}_*\sigma_E\ra,$$
where $c_1$ denotes the first Chern class and $\pi_{S^1}:F_{S^1}\to
BS^1$ is the projection.

So fix a homology class $\beta\in H_2^{S^1}(F;\ZZ)$. Let
$E\to\Sigma$ be the principal $S^1$-bundle of degree
$\la c_1(ES^1),{\pi_{S^1}}_*\beta\ra$.
Define as in Section \ref{theequations} $\AAA$ to be the set of
connections on $E$ and $\SSS$ the set of smooth sections 
of $\FFF=E\times_{S^1}F$. 
Let $$\bM(\beta,c)=\bM_I(\beta,c)=\{ (A,\phi)\in\AAA\times\SSS\mid
{\phi_E}_*[\Sigma]=\beta\text{ and $(A,\phi)$ satisfies
  (\ref{equs})}\},$$
where $\phi_E:\Sigma\to F_{S^1}$ is the map constructed in Lemma
\ref{bijeccio}. 
By gauge equivariance of the equations and of the condition
${\phi_E}_*[\Sigma]=\beta$ (by Lemma \ref{bijeccio}, 
if ${\phi_E}_*[\Sigma]=\beta$ and $g\in\GGG$
then $(g^*\phi_E)_*[\Sigma]=\beta$), 
it makes sense to define
$$\MMM(\beta,c)=\MMM_I(\beta,c)=\bM_I(\beta,c)/\GGG.$$ 
We will only use the notations 
$\bM_I(\beta,c)$ and $\MMM_I(\beta,c)$ when the complex structure $I$
is not clear from the context.

\subsection{The choice of $c$}
Assume that the action of $S^1$ on $F$ is semi-free (i.e., $S^1$ acts freely
on the complementary $F\setminus F^f$ of the fixed point sets).
Let $F_1,\dots,F_r\subset F$ be the connected components of the fixed
point set $F^f$. It follows from the properties of
the moment map that for any $1\leq j\leq r$ the restriction $\mu|_{F_j}$ 
takes a constant value, say $c_j\in\imag\RR$.
Let
$$\CCC=\{ c_j-2\pi\imag \deg(E)/\Vol(\Sigma)\mid 1\leq j\leq r\}.$$

\begin{lemma}
If $c\in\imag\RR\setminus\CCC$, then the action of $\GGG$ on 
$\bM(\beta,c)$ is free.
\label{elsbonsvalors}
\end{lemma}
\begin{pf}
Take any $(A,\phi)\in\AAA\times\SSS$, and assume that there is a
nontrivial gauge transformation which fixes $(A,\phi)$. Since the
stabiliser of any connection is the set of constant gauge
transformations, we deduce that $(A,\phi)$ is fixed by a nontrivial
constant gauge. 
Because the action of $S^1$ on $F$ is semi-free, this implies that
$\phi(\Sigma)\subset E\times_{S^1}F^f\subset\FFF$, and 
since $\Sigma$ is connected, we deduce that 
\begin{equation}
\phi(\Sigma)\subset E\times_{S^1}F_j
\label{inclosfix}
\end{equation}
for some $1\leq j\leq r$.
Now assume that $(A,\phi)$ satisfies
$$\Lambda F_A+\mu(\phi)=c$$
for some $c\in\imag\RR$. Integrating this equality over $\Sigma$,
dividing by $\Vol(\Sigma)$ and using Chern--Weil theory to
write $\frac{\imag}{2\pi}\int_{\Sigma} F_A=\deg(E)$, we deduce from
(\ref{inclosfix}) that 
$$c=c_j-2\pi\imag \deg(E)/\Vol(\Sigma)\in\CCC.$$
\end{pf}

\subsection{The universal bundle}
\label{univbundle}
Let 
$$\bEE=\pi_{\Sigma}^*E=
\AAA\times\SSS\times E\to \AAA\times\SSS\times\Sigma,$$
where $\pi_{\Sigma}:\AAA\times\SSS\times\Sigma\to\Sigma$ is the projection.
Consider the section of the associated bundle
$$\begin{array}{rcl}
\obPhi:\AAA\times\SSS\times\Sigma & \to & \bEE\times_{S^1}F
=\AAA\times\SSS\times\FFF \\ 
(A,\phi,x) & \mapsto & (A,\phi,\phi(x)).
\end{array}$$
Take on $\AAA\times\SSS\times\Sigma$ the diagonal action of $\GGG$,
acting trivially on $\Sigma$. There is a canonical lift of
this action to $\bEE$ (resp. to $\bEE\times_{S^1}F$) which is the
diagonal action on $\AAA\times\SSS\times E$ (resp. 
$\AAA\times\SSS\times\FFF$), and the section $\obPhi$ is $\GGG$-equivariant.

Suppose that $c\in\imag\RR\setminus\CCC$. By Lemma \ref{elsbonsvalors}
the action of $\GGG$ on $\bM(\beta,c)\times\Sigma$ is free. This implies
that the restriction of $\bEE$ to $\bM(\beta,c)\times\Sigma$ descends to
give a bundle $$\EE\to\MMM(\beta,c)\times\Sigma.$$
We call $\EE$ the universal bundle. By equivariance the
section $\obPhi$ also descends to give a section $\bPhi$ of the
bundle $\EE\times_{S^1}F\to\MMM(\beta,c)\times\Sigma$. 
Arguing like in Lemma \ref{bijeccio} we obtain out of $\EE$ and $\bPhi$
a canonical (homotopy class of)
map $$\ev_{\beta,c}:\MMM(\beta,c)\times\Sigma\to F_{S^1},$$
which we call the equivariant evaluation map.

\subsection{The invariants}
We now give a heuristic definition of the invariants.
Suppose that $c\in\imag\RR\setminus\CCC$, so that by the preceeding
section we have an equivariant evaluation map
$\ev_{\beta,c}:\MMM(\beta,c)\times\Sigma\to F_{S^1}$.
Suppose that there is a {\it canonical} fundamental
class $$[\MMM(\beta,c)]\in H_*(\MMM(\beta,c);\ZZ).$$ 
Let $\alpha_1,\dots,\alpha_p\in H^*_{S^1}(F;\ZZ)$ be equivariant cohomology
classes and let $\gamma\in H^*(\AAA/\GGG;\ZZ)$. Let
$\nu:\MMM\to\AAA/\GGG$ be the map induced by the projection
$\AAA\times\SSS\to\AAA$.
We define the Hamiltonian Gromov--Witten invariant
$\Psi_{\Sigma,\beta,c}$ evaluated at $\alpha_1,\dots,\alpha_p$ and
$\gamma$ to be
\begin{equation}
\Psi_{\Sigma,\beta,c}(\alpha_1,\dots,\alpha_p;\gamma)
=\la\ev_{\beta,c}^*\alpha_1/[\Sigma]\cup\dots\cup
\ev_{\beta,c}^*\alpha_p/[\Sigma]\cup\nu^*\gamma,
[\MMM(\beta,c)]\ra\in\ZZ.
\label{definv}
\end{equation}
Of course, the existence of a canonical fundamental class is far from
obvious, since in general $\MMM(\beta,c)$ is neither compact nor a
smooth manifold. However, perturbing slightly the equations
(\ref{equs}) we obtain a smooth set of gauge orbits of solutions,
which can be compactified in a natural way. And although we will not
construct any fundamental homology class of the compactification, in
some situations we will manage to give a rigorous definition of the
invariants sketched above using the technique of pseudo-cycles.

\subsection{Assumptions on $F$ and the action of $S^1$}
In order to give a rigorous definition of the invariants 
we will assume that $F$ and the action of $S^1$ satisfy several
technical conditions, which we now list.
\begin{enumerate}
\item The action of $S^1$ on $F$ is semi-free; 
\item the manifold $F$ is monotone; this means that there is a
 real number $\lambda > 0$ such that $[\omega]=\lambda c_1(TF)$;
\item the connected components $F_1,\dots,F_r$ of the fixed point set
  $F^f$ are also monotone, i.e., there exists positive real numbers
  $\lambda_1,\dots,\lambda_r>0$ such that for any $1\leq k\leq r$
  $$[\omega_k]=\lambda_k c_1(TF_k),$$ where $\omega_k\in\Omega^2(F_k)$
  is the restriction of $\omega$; furthermore, for any  $1\leq k\leq
  r$ we also have $$\lambda_k\leq\lambda;$$
\item we have for any $1\leq k\leq r$ a bound
\begin{equation}
\codim_{\RR}F_k\leq 6.
\label{cond3}
\end{equation}
\end{enumerate}
For example, any monotone symplectic manifold of dimension $\leq 6$
with a Hamiltonian action with isolated fixed points satisfies the
above conditions.

We will use the following consequences of the above assumptions. Let
$s:\CP^1\to F$ be a $I$-holomorphic map. Then
\begin{equation}
0<\la c_1(s^*TF),[\CP^1]\ra.
\label{cond1}
\end{equation}
Furthermore, if $s(\CP^1)\subset F_k$, then 
\begin{equation}
\la c_1(s^*TF_k),[\CP^1]\ra\leq\la c_1(s^*TF),[\CP^1]\ra.
\label{cond2}
\end{equation}
Both (\ref{cond1}) and (\ref{cond2}) follow from our assumptions, 
thanks to the formula
$$\int_{\CP^1}|ds|^2=\int_{\CP^1}s^*\omega,$$
which is valid because $s$ is $I$-holomorphic.
\label{conditions}

\noindent{\bf Remarks.}
Some comments on the technical conditions are now in order.

\noindent{\bf 1.} The condition that the action is semi-free is used
in two places. First to assure that there is universial bundle
$\EE\to\MMM\times\Sigma$, and then to prove that a generic perturbation
of the equation gives a smooth moduli of solutions. 
If the action is not semi-free, one might encounter two different (but
related) problems. First, the moduli of solutions $\MMM$ could be
forced to be an orbifold, and not a manifold, and there would
not exist a universal bundle, but an orbibundle. Second, the kind of
perturbations which we use will not suffice to get smoothness, and
multivalued perturbations will have to be used instead (see \cite{Sa}).
Finally, we recall that, although the condition of being semi-free is
rather restrictive, it has also a very natural geometrical 
meaning: namely, it is equivalent to the condition of all 
Marsden--Weinstein quotients at regular values of $\mu$ being smooth.

\noindent{\bf 2.} Conditions (\ref{cond1}), (\ref{cond2}) and
(\ref{cond3}) are of technical nature, and
using more sophisticated techniques such as as virtual moduli cycles one
might presumably get rid of them.

\noindent{\bf 3.} It seems natural that one could define invariants of
Hamiltonian actions of arbitrary compact Lie groups using the same
ideas given as in this paper. However, giving a rigorous definition of
them (at least using the techniques which we deal with here) seems to
involve a considerable amount of extra work. In \cite{Mu2}, parts of
the programme developped in this paper are also worked out for
arbitary compact Lie groups (for example, the compactification of the
moduli of solutions). 
On the other hand, note that if we replace $S^1$ by tre trivial group
we recover the usual Gromov--Witten invariants.

\subsection{Kaehler situation}
The second equation in (\ref{equs}) was studied in \cite{Mu1} under the
assumption that the complex structure $I$ is integrable (i.e., when it
gives a Kaehler structure to $F$). In this case the action of $S^1$
extends to a holomorphic action of $\CC^*$. This allows to extend the
action of $\GGG$ on $\SSS$ to an action of
$\GGG^{\CC}=\Map(\Sigma,\CC^*)$. On the other hand, there is a natural
extension of the action of $\GGG$ on $\AAA$ to an action of
$\GGG^{\CC}$. The theorem proved in \cite{Mu1} describes which orbits
of the action of $\GGG^{\CC}$ on $\AAA\times\SSS$
contain solutions to the second equation
in (\ref{equs}). When the Kaehler manifold $(F,\omega,I)$ is
projective, the description might be given using concepts coming from
GIT. This description was used in \cite{Mu2} to compute some nonzero
Hamiltonian Gromov--Witten invariants of $S^2$ with the action of
$S^1$ given by rotations.

\subsection{Symplectic interpretation}
As most gauge theoretical equations, the second equation in
(\ref{equs}) admits an infinite dimensional symplectic interpretation,
which we now briefly explain (see \cite{Mu1} for more details). The
space $\AAA\times\SSS$ admits a natural symplectic structure (which is
the sum of the symplectic structure in $\AAA$ defined by M. Atiyah and
R. Bott in \cite{AtBo} and a symplectic structure on $\SSS$ obtained
using $\omega$). Then the action of $\GGG$ on $\AAA\times\SSS$
is Hamiltonian, and
$$\nu(A,\phi):=\Lambda
F_A+\mu(\phi)\in\Omega^0(\imag\RR^*)\subset\Omega^0(\imag\RR)^*$$ is a
moment map for this action. On the other hand, the set
$$\DDD_{\beta}=
\{(A,\phi)\mid \ov{\partial}_A\phi=0,\ {\phi_E}_*[\Sigma]=\beta\}$$
is $\GGG$-invariant. So $\MMM(\beta,c)$ is the symplectic quotient of
$\DDD_{\beta}$ at the central element $c\in\Omega^0(\imag\RR)^*$ (by a slight
abuse of notation here $c$ denotes the constant section with value
$c$). Hence, the smooth locus of $\MMM(\beta,c)$ has a natural
symplectic structure.
On the other hand, using the action of $\GGG$ on the restriction
$\bE|_{\DDD_{\beta}\times\Sigma}$, we get a $S^1$-principal bundle on
the Borel construction $$\EE_{\GGG}\to
E\GGG\times_{\GGG}(\DDD_{\beta}\times\Sigma)$$
and, proceeding exactly as in the construction of $\bPhi$  
in Section \ref{univbundle}, we can construct a section
$\bPhi_{\GGG}\in\Gamma(\EE_{\GGG}\times_{S^1}F).$
Using this data we get, by means of Lemma \ref{bijeccio}, a map
$$\ev_{\GGG}:E\GGG\times_{\GGG}(\DDD_\beta\times\Sigma)\to F_{S^1}.$$
Now, it is straightforward to prove that
$$\ev_{\beta,c}^*=\kappa_c\circ\ev_{\GGG}^*,$$
where $\kappa_c:H_{\GGG}^*(\DDD_\beta\times\Sigma)\to
H^*(\MMM(\beta,c)\times\Sigma)$ is the Kirwan map.

\section{Smoothness of the moduli space}
\label{s:smoothness}
\subsection{Sobolev completions}
To be able to deal with Banach manifold techniques (such as 
Sard--Smale theorem) we will work with the
completion of the space $\AAA\times\SSS$ with respect to suitable 
Sobolev norms. Take any $p>2$. Let $A_0\in\AAA$ be any smooth
connection, and define
$\AAA_{L^p_1}=A_0+L^p_1(T^*\Sigma\otimes\imag\RR)$
(if $E\to\Sigma$ is a vector bundle with a metric $|\cdot|$ and a connection
$\nabla$, we write $L^p_1(E)$ for the completion of 
$\Omega^0(\Sigma;E)$ with respect to
the norm $\|\sigma\|_{L^p_1}=\|\sigma\|_{L^p}+\|\nabla\sigma\|_{L^p}$;
and to define $L^p_1(T^*\Sigma\otimes\imag\RR)$ we use any connection
on $T^*\Sigma$).
It is easy to check that $\AAA_{L^p_1}$ is independent of $A_0$.
Take any smooth 
embedding $j:\FFF\to\RR^N$; any $\phi\in\SSS$ can be seen as a map
$\phi:\Sigma\to\FFF$ and we set, for any $\phi,\phi'\in\SSS$,
$$d_p(\phi,\phi')=\|j(\phi)-j(\phi')\|_{L^p_1}=
\sum_{1\leq k\leq N} \|e_k(j(\phi))-e_k(j(\phi'))\|_{L^p_1},$$ 
where $e_k:\RR^N\to\RR$ is the $k$-th coordinate.
Let $\SSS_{L^p_1}$ be the completion of $\SSS$ with respect to the
metric $d_p$.
Since $\dim_{\RR}\Sigma=2$, the Sobolev theorems tell us that
there is an inclusion $L^p_1(\Sigma)\subset C^0(\Sigma)$. 
Hence, the elements
of $\SSS_{L^p_1}$ are continuous sections. Furthermore, $\SSS_{L^p_1}$
is independent of the embedding $j$.
Finally, we complete $\GGG=\Map(\Sigma,S^1)$ using $L^p_2$ norm.
Then $\GGG_{L^p_2}$ is a Banach Lie group, which acts smoothly on
$\AAA_{L^p_1}$ and $\SSS_{L^p_1}$.

Let $\SSS^*=\{\phi\in\SSS\mid \phi(\Sigma)\nsubset
E\times_{S^1}F^f\}$, and let $\SSS^*_{L^p_1}$
be the closure of $\SSS^*$ in $\SSS_{L^p_1}$.
One can give a natural Banach manifold structure to the set
$$\BBB=(\AAA_{L^p_1}\times\SSS^*_{L^p_1})/\GGG_{L^p_2}.$$
There is also a natural Banach bundle structure on 
$\bW\to\AAA_{L^p_1}\times\SSS^*_{L^p_1}$, 
where the fibre over $(A,\phi)$ is
$$\bW_{(A,\phi)}=L^p(T^{0,1}\Sigma\otimes\phi^*T\FFF_v)
\oplus L^p(\imag\RR).$$   
The action of $\GGG_{L^p_2}$ lifts to a smooth action on 
$\bW$ and there is a quotient Banach bundle
$\WWW\to\BBB$.
Finally, equations (\ref{equs}) provide a smooth section
$\Psi:\BBB\to\WWW$. 

In the sequel we will omit the Sobolev subscripts in $\AAA$, $\SSS$
and $\GGG$, and Sobolev completions will be implicity assumed unless
otherwise stated.

\subsection{Virtual dimension of the moduli}
The section $\Psi$ is Fredholm, that is, its covariant derivative
at any $[A,\phi]\in\Psi^{-1}(0)\subset\BBB$ (here the brackets
denote gauge equivalence class) is a Fredholm operator
$D_{[A,\phi]}\Psi:T_{[A,\phi]}\BBB\to\WWW_{[A,\phi]}$.
Its index (we will be talking on real indices unless 
otherwise specified)
is equal to minus the index of the deformation complex 
\begin{equation}
C^0_{A,\phi}\stackrel{d_0}{\longrightarrow}
C^1_{A,\phi}\stackrel{d_1}{\longrightarrow}
C^2_{A,\phi},
\label{complexdefo}
\end{equation}
where $C^0_{A,\phi}=L^p_2(\imag\RR)=\Lie\GGG$,
$C^1_{A,\phi}=L^p_1(\phi^*T\FFF_v)\oplus L^p_1(T^*\Sigma\otimes\imag\RR)=
T_{(A,\phi)}\AAA\times\SSS$ and
$C^2_{A,\phi}=L^p(T^{0,1}\Sigma\otimes\phi^*T\FFF_v)
\oplus L^p(\imag\RR)$; $d_0$ is given by the infinitesmial action 
of $\GGG$ on $\AAA\times\SSS$ and $d_1$ is the linearisation of
equations (\ref{equs}).
Modulo compact operators, the complex (\ref{complexdefo}) splits
as the sum of the de Rham complex
$$L^p_2(\imag\RR)\stackrel{d}{\longrightarrow}
L^p_1(T^*\Sigma\otimes\imag\RR)\stackrel{d'}{\longrightarrow}
L^p(\Lambda^2T^*\Sigma\otimes\imag\RR)\simeq L^p(\imag\RR)$$
(the isomorphism being given by $\Lambda$), whose index is 
$$\Ind d'+d^*=2-2g,$$ 
plus a (shifted one unity) Dolbeault complex
$$L^p_1(\phi^*T\FFF_v)
\stackrel{\ov{\partial}_{\phi}}{\longrightarrow}
L^p(T^{0,1}\Sigma\otimes\phi^*T\FFF_v),$$
whose index is, by Riemann--Roch, equal to 
\begin{align*}
\Ind \ov{\partial}_{\phi}&=2\la c_1(\phi^*T\FFF_v),[\Sigma]\ra+2n(1-g) \\
&= 2\la \phi^*c_1(T\FFF_v),[\Sigma]\ra+2n(1-g).
\end{align*}  
Consequently, we have
$\Ind D_{[A,\phi]}\Psi=
2\la \phi^*c_1(T\FFF_v),[\Sigma]\ra+2(n-1)(1-g).$
Observe that this is a function of
$\beta={\phi_E}_*[\Sigma]$. More precisely, we may write
\begin{equation}
\Ind D_{[A,\phi]}\Psi=2\la c_1^{S^1}(TF),\beta\ra+2(n-1)(1-g),
\label{dimvirtual}
\end{equation}
where $c_1^{S^1}$ denotes the first equivariant Chern class.
As usual, we call this number the real virtual dimension of our moduli
space.

\subsection{Perturbing the equations}
Recall that $\pi:\FFF\to\Sigma$ denotes the projection.
Let $\Hom^{0,1}(\pi^*T\Sigma,T\FFF_v)$ be the space of antiholomorphic 
smooth vector bundle maps from $\pi^*T\Sigma$ to $T\FFF_v$. 
This space has an action of $S^1$ induced 
by the action on $\FFF$. Let
$$\PPP=\Hom^{0,1}(\pi^*T\Sigma,T\FFF_v)^{S^1}
\oplus\Omega^0(\Sigma;\imag\RR)$$ 
(the superscript $S^1$ denotes the subspace of invariant
elements). For any $\rho=(\rho_1,\rho_2)\in\PPP$, consider the
following perturbed equations
\begin{equation}
\left\{\begin{array}{l}
\ov{\partial}_A\phi=\rho_1\\
\Lambda F_A+\mu(\phi)=c+\rho_2
\end{array}\right.
\label{equs2}
\end{equation}
and define the perturbed set of solutions to be
$$\bM_\rho(\beta,c)=\{(A,\phi)\in\AAA\times\SSS\mid
{\phi_E}_*[\Sigma]=\beta\text{ and $(A,\phi)$ satisfies (\ref{equs2})}\}.$$
Since the elements of $\PPP$ are $S^1$-invariant, the set
$\bM_\rho(\beta,c)$ is gauge invariant, so we may define the perturbed
moduli space to be 
$$\MMM_{\rho}(\beta,c)=\bM_\rho(\beta,c)/\GGG.$$
In order to avoid having solutions of the perturbed equations which
are entirely contained in the fixed point set, we assume that
$c\in\imag\RR\setminus\CC$ and we define
$$\PPP_c=\{(\rho_1,\rho_2)\in\PPP\mid |\rho_2|<d(c,\CCC)\},$$
where $d(c,\CCC)$ denotes the distance from $c$ to the set
$\CCC$. Then we may prove, exactly like in Lemma \ref{elsbonsvalors}:
\begin{lemma}
If $\rho\in\PPP_c$, then the action of $\GGG$ on 
$\bM_\rho(\beta,c)$ is free.
\label{elsbonsvalors2}
\end{lemma}
The following Theorem justifies the use of the perturbed equations.

\begin{theorem}
Assume that $c\in\imag\RR\setminus\CCC$. There is a subset
$\PPP_c^{\reg}\subset\PPP_c$ of Baire of the second category 
(with respect to the $C^{\infty}$ topology on $\PPP_c$) such that
for any $\rho\in\PPP_c^{\reg}$ the perturbed moduli space 
$\MMM_{\rho}(\beta,c)$ is a smooth manifold of dimension equal to
(\ref{dimvirtual}) with a natural orientation. 
Furthermore, for any pair of perturbations
$\rho,\rho'\in\PPP_c$ there is a path $P:[0,1]\to\PPP_c$ with
$P(0)=\rho$, $P(1)=\rho'$ and such that
$$\MMM_P(\beta,c)=\bigcup_{t\in[0,1]}\MMM_{P(t)}(\beta,c)$$
is a smooth oriented cobordism between $\MMM_{\rho}(\beta,c)$, and
$\MMM_{\rho'}(\beta,c)$.
\label{modulipert}
\end{theorem}
\begin{pf}
The techniques needed to prove this result are rather standard, so 
we will be very sketchy. More details are given in \cite{Mu2}, Theorem
3.4.4, and see also \cite{FrUh, McDS1} for analogous results.

Take a big enough (to be specified later)
integer $l>0$, and consider the completion $\PPP_c^l$ of
$\PPP_c$ with respect to the $C^l$ norm. Let
$$\MMM_{\PPP}(\beta,c)=\{(\rho,A,\phi)\in\PPP_c^l\times\AAA\times\SSS
\mid{\phi_E}_*[\Sigma]=\beta\text{ and $(A,\phi,rho)$ 
satisfies (\ref{equs2})}\}/\GGG.$$
One first proves that this is a Banach manifold.
Indeed, consider the following section of the pullback bundle 
$\pi_{\BBB}^*\WWW\to\PPP_c^l\times\BBB$ (where $\pi_{\BBB}$
denotes the projection):
$$\Psi_{\PPP}((\rho_1,\rho_2),[A,\phi])=(\ov{\partial}_A\phi-\rho_1,\Lambda
F_A+\mu(\phi)-c-\rho_2)$$ 
(here $[A,\phi]$ denotes de gauge equivalence
class of $(A,\phi)\in\AAA\times\SSS$). Then 
$\MMM_{\PPP}(\beta,c)=\Psi_{\PPP}^{-1}(0)$. So to prove that 
$\bM_{\PPP}(\beta,c)$ is a Banach manifold it suffices to prove that
if $z=((\rho_1,\rho_2),[A,\phi])\in\Psi_{\PPP}^{-1}(0)$ then
the differential $D\Psi_{\PPP}(z):T(\PPP_c^l\times\BBB)\to
(\pi_{\BBB}^*\WWW)_z$ is onto. 
By ellipticity, we know that the image of $D\Psi_{\PPP}(z)$ is
closed. So if $D\Psi_{\PPP}(z)$ were not onto, there would be some
nonzero $\eta\in L^q(T^{0,1}\Sigma\otimes\phi^*T\FFF_v)\oplus L^q(\imag\RR)$,
where $p^{-1}+q^{-1}=1$, and such that
$\la\eta,D\Psi_{\PPP}(z)\xi\ra=0$ for any $\xi\in
T(\PPP_c^l\times\BBB)$. Now, it is easy to prove, using that
$\phi\in\SSS^*$ (here we use Lemma \ref{elsbonsvalors}), 
that one can find some $\xi\in T\PPP_c^l\subset
T(\PPP_c^l\times\BBB)$ which does not satisfy that equality. This
proves that $\bM_{\PPP}(\beta,c)$ is a Banach manifold.

Consider the projection 
$\pi_{\PPP}:\MMM_{\PPP}(\beta,c)\to\PPP^l$. This is a map whose
differential is everywhere
Fredholm, and its index $\Ind(D\pi_{\PPP})$
is equal, at any point, to the virtual dimension (\ref{dimvirtual}).
Now, provided $l>2+\Ind(D\pi_{\PPP})$,
we may apply Sard--Smale theorem to deduce that there is a set of
Baire of the second category $\PPP_c^{\reg,l}\subset\PPP_c^l$
of regular values of $\pi_{\PPP}$. And $\pi_{\PPP}^{-1}(\rho)=
\MMM_{\rho}(\beta,c)$ by definition. Finally, a trick of Taubes
(see p. 36 in \cite{McDS1})
allows to deduce from this result that there is also a subset
$\PPP_c^{\reg}\subset\PPP_c$ of Baire of the second category
(with respect to the $C^{\infty}$ topology on $\PPP_c$) of regular
values of $\pi_{\PPP}$. 

The result on cobordisms is proved in the same way. 

Finally, the orientability follows from identifying the tangent
space at $[A,\phi]\in\MMM_{\rho}(\beta,c)$ with the first cohomology
group of the complex (\ref{complexdefo}). By Hodge theory this group
can be identified with the kernel of the elliptic operator $d_0^*+d_1$
and, when $\rho\in\PPP_c^{\reg}$, this kernel carries a natural orientation
because $d_0^*+d_1$ has the same symbol as a
Cauchy--Riemann operator plus the Hodge operator $d+d^*$ acting on
1-forms.
\end{pf}

\subsection{Regularity}
\label{sb:regularity}
The following theorem proves that the peturbed moduli spaces 
$\MMM_{\rho}(\beta,c)$ which we get by taking smooth perturbations
$\rho\in\PPP$ are independent of the chosen Sobolev completion.

\begin{theorem}
Let $\rho\in\PPP$, and let
$(A,\phi)\in\AAA\times\SSS$ be a solution to the perturbed
equations (\ref{equs2}). There is a gauge transformation $g\in\GGG$
such that $g^*A$ and $g^*\phi$ are both smooth.
\label{regularitat}
\end{theorem}
\begin{pf}
Take $g\in\GGG$ such that $g^*A$ is in Coulomb gauge with respect to
the smooth connection $A_0\in\AAA$, i.e., such that
$d^*(A_0-g^*A)=0$, and define $A_s=g^*A$ and $\phi_s=g^*\phi$.
We will prove that $A_s$ and $\phi_s$ are smooth. To begin with we
know that the $L^p_1$ norms of $A_0-A_s$ and of $\phi_s$ are bounded.
Let $\omega_{\Sigma}$ be the volume form of $\Sigma$. The second
equation in (\ref{equs2}) may be written
$$d(A_0-A_s)=-F_{A_0}+\omega_{\Sigma}(c-\mu(\phi)+\rho_2).$$
Since $\mu$ is smooth and $\phi$ has bounded $L^p_1$ norm, we deduce
that the $L^p_1$ norm of $\mu(\phi)$ is also bounded. So the above
equation gives an $L^p_1$ bound to $d(A_0-A_s)$. This bound, combined
with $d^*(A_0-A_s)=0$ implies, by the ellipticity of $d+d^*$, an
$L^p_2$ norm on $A_0-A_s$. From this we obtain a bound on the $L^p_2$
norm of the complex structure $I(A_s)$ constructed in Subsection
\ref{corbeshol}. Now, standard results on regularity of (perturbed)
holomorphicity equation for curves (see for example Theorem B.3.4 in
\cite{McDS1}), allow to deduce from this $L^p_2$ bounds on the section
$\phi$ (which is a $I(A_s)$ holomorphic map from $\Sigma$ to $\FFF$, by
Lemma \ref{lcorbeshol}). So we have passed from $L^p_1$ bounds on $A_s$
and $\phi_s$ to $L^p_2$ bounds on both. This step can be repeated to
obtain $L^p_k$ bounds for any $k$. And this implies, by the Sobolev
theorems, that both $A_s$ and $\phi_s$ are smooth.  
\end{pf}

What this theorem proves is, strictly speaking, that 
$\MMM_{\rho}(\beta,c)$ is independent of $p$ as a set. To prove also
that, when $\rho\in\PPP_c^{\reg}$, 
the structure of $\MMM_{\rho}(\beta,c)$
as a differentiable manifold is intrinsic, one may
use standard elliptic theory applied to the Kuranishi models which
describe locally the moduli spaces (see Section 3.3 in \cite{Mu2}).

\section{Compactification of the moduli space}
\label{s:compactification}

\subsection{The Yang--Mills--Higgs functional}
Define the Yang--Mills--Higgs functional $\YMH_c:\AAA\times\SSS\to\RR$
as
$$\YMH_c(A,\phi)=\|F_A\|_{L^2}^2+\|d_A\phi\|_{L^2}^2+
\|\mu(\phi)-c\|_{L^2}^2$$ for any $(A,\phi)\in\AAA\times\SSS$.
The following is Lemma 7.9 in \cite{Mu1} (note that in \cite{Mu1} we
assume the manifold $F$ to be Kaehler; however, the results on the
functional $\YMH_c$ remain valid when the complex structure on $F$ is
not integrable).

\begin{lemma}
For any $(A,\phi)\in\AAA\times\SSS$ we have
$$\YMH_c(A,\phi)=\|\Lambda F_A+\mu(\phi)-c\|_{L^2}^2+
2\|\ov{\partial}_A\phi\|_{L^2}^2+\int_{\Sigma}\la F_A,c\ra+
\la \phi^*[\omega_{\FFF}],\Sigma\ra,$$
where $[\omega_{\FFF}]\in H^2(\FFF;\RR)$ 
is the cohomology class represented
by the coupling form of the symplectic fibration $\FFF\to\Sigma$
and the connection $A$
(this cohomology class does not depend on $A$, see \cite{GuLeS}).
\end{lemma}

\begin{corollary}
Given $\beta\in H_2^{S^1}(F;\ZZ)$, $c\in\imag\RR$ and $\rho\in\PPP$, 
there is a constant $C(\beta,c,\rho)>0$ such that for any 
$(A,\phi)\in\bM_{\rho}(\beta,c)$
we have $\|F_A\|_{L^2}<C(\beta,c,\rho)$ and 
$\|d_A\phi\|_{L^2}<C(\beta,c,\rho)$.
\label{energybound}
\end{corollary}
\begin{pf}
Indeed, the term 
$\int_{\Sigma}\la F_A,c\ra+\la \phi^*[\omega_{\FFF}],\Sigma\ra$ only
depends on $\beta$ and $c$.
\end{pf}
In fact, the bound for $F_A$ is obvious when the manifold
$F$ is compact, since then $\mu$ is bounded; 
although we will not need it, we mention that 
the result remains valid even when $F$ is not compact.

\subsection{Compactness}
Let $\rho\in\PPP$ be any perturbation.
                                   
\begin{definition} A cusp $\rho$-THC is the following set of data.
\begin{enumerate}
\item A compact connected singular curve $\Sigma^{\cusp}$ with only nodal 
singularities, of the form
$\Sigma^{\cusp}=\Sigma_0\cup\Sigma_1\cup\dots\cup\Sigma_K$, 
where $\Sigma_0=\Sigma$ is called the principal component, 
and where the other components are rational curves $\CP^1$
and are called bubbles; furthermore, two different components
$\Sigma_i$ and $\Sigma_j$ meet at most at one point. 
\item A $S^1$-principal bundle $E\to\Sigma_0$, a connection $A$ on $E$,
a section $\phi_0:\Sigma_0\to \FFF=E\times_{S^1}F$ and an element
$c\in\imag\RR$ satisfying the equations
$$\left\{\begin{array}{l}
\ov{\partial}_A\phi_0=\rho_1 \\
\Lambda F_A+\mu(\phi_0)=c+\rho_2. \end{array}\right.$$
\item For any $k\neq 0$, a holomorphic map $\phi_k:\Sigma_k\to\FFF$ whose 
image is inside a unique fibre $\FFF_{x_k}$ of $\FFF\to\Sigma$ (note $\phi_k$
is holomorphic with respect to the complex structure on $F$). The maps
$\phi_0,\phi_1,\dots,\phi_K$ are required to glue together to give
a map $\phi:X^{\cusp}\to\FFF$. 
\end{enumerate}
\label{defcuspTHC}
\end{definition}
We denote cusp $\rho$-THCs with tuples of the form 
$(\Sigma^{\cusp},E,A,\phi,c)$.
Let $\iota:H_*(\FFF_x;\ZZ)\to H_*^{S^1}(F;\ZZ)$ be the natural map
from the homology of any fibre of $\FFF$ to the equivariant homology
of $F$ (this map is well defined because, since $S^1$ is connected,
the action of $S^1$ on $H_*(F;\ZZ)$ is trivial).
We will say that the $\rho$-THC $(\Sigma^{\cusp},E,A,\phi,c)$
represents the class
$$({\phi_0}_E)_*[\Sigma_0]+\sum_{k=1}^K\iota_*{\phi_k}_*[\Sigma_k]\in
H_2^{S^1}(F;\ZZ).$$ 

\begin{theorem}
Let $\beta\in H_2^{S^1}(F;\ZZ)$ 
be any homology class. 
Consider a sequence of gauge equivalence classes
$[A_j,\phi_j]\in\MMM_{\rho}(\beta,c_j)$, where $j\geq 1$, and 
where $\{c_j\}\subset\imag\RR$ is a bounded set.
After passing to a subsequence, there exists a cusp $\rho$-THC 
$(\Sigma^{\cusp},E,A,\phi,c)$ and gauge transformations $g_j\in\GGG$ 
such that if $(A_j',\phi_j')=g_j(A_j,\phi_j)$ we have
\begin{enumerate}
\item $c_j\to c$;
\item $A_j'\to A$ in $C^{\infty}$;
\item the $\rho$-THC $(\Sigma^{\cusp},E,A,\phi,c)$ represents the
class $\beta$;  
\item the images $\phi_j'(X)\subset\FFF$ converge pointwise to
$\phi(X^{\cusp})$, that is, for any sequence $x_j\in X_j$
there exists $x\in X^{\cusp}$ such that $\phi_j'(x_j)\to \phi(x)$.
\end{enumerate}
Furthermore, the curve $\Sigma^{\cusp}$ is a tree, in the sense that
the graph with one point for each irreducible component of
$\Sigma^{\cusp}$ and with an edge joining two points exactly when
these correspond to components which intersect, is a tree. And,
finally, the limit $\rho$-THC curve $(\Sigma^{\cusp},E,A,\phi,c)$
is stable, which means that if the map $\phi_k$ for $1\leq k\leq K$ is
constant, then the bubble $\Sigma_k$ meets at least three other
irreducible components of $\Sigma^{\cusp}$.
\label{compactificacio}
\end{theorem}
\begin{pf}
We first take a subsequence of $(A_j,\phi_j,c_j)$ (and denote it with the
same symbol) such that $c_j\to c\in\imag\RR$.
Just as in the proof of Theorem \ref{regularitat}, we take gauge
transformations which put the connections in Coulomb gauge and then
use the existing compactness theorem for pseudoholomorphic curves.
So let $g_j\in\GGG$ such that $A_j'=g^*A_j$ satisfies
$d^*(A'_j-A_0)=0$, and let $\phi_j'=g^*\phi_j$. By gauge invariance,
$(A_j',\phi_j')$ satisfies
$$\left\{\begin{array}{l}
\ov{\partial}_A\phi_0=\rho_1 \\
\Lambda F_A+\mu(\phi_0)=c+\rho_2. \end{array}\right.$$
Now, the second equation combined with the Coulomb condition gives a
uniform bound $\|A_j'-A_0\|_{L^p_1}<C$, by ellipticity of $d+d^*$.
Using Rellich theorem on the compactness of the embedding $L^p_1\to
C^0$, we deduce that there is a connection $A$ satisfying
$\|A-A_0\|_{L^p_1}<\infty$ such that, after restricting to a 
subsequence (which we denote also by $\{A_j',\phi_j'\}$), 
$\{A_j'\}$ converges uniformly in $C^0$ to $A$.
This implies that the complex structures $I_j=I(A_j)$ converge 
in $C^0$ to $I=I(A)$. 
On the other hand, $\phi'_j:\Sigma\to\FFF$ is a perturbed 
$I_j$ holomorphic map for any $j$ (by Lemma \ref{lcorbeshol}). 
(These maps can be regarded also as
holomorphic maps to $\Sigma\times\FFF$, giving to this space a complex
structure of the form $\left(\begin{array}{cc}
I_{\Sigma} & 0 \\ \rho_1 & I_j\end{array}\right)$, see for example
\cite{Gr,Ru1}.) By Corollary \ref{energybound} there
are bounds on $\|d_{A_j'}\phi_j'\|_{L^2}^2$. But if we give to $\FFF$
the Riemannian metric $g_j$ obtained by summing the metric
$\omega(\cdot,I\cdot)$ on $F$ and the metric on $\Sigma$ by means of
the splitting of $T\FFF$ given by the connection $A_j'$, then we have
$$\|d\phi_j'\|_{L^2,g_j}^2=\|d_{A_j'}\phi_j'\|_{L^2}^2+\Vol(\Sigma),$$
so the energy of the maps $\phi_j'$ is bounded. On the other hand, 
the metrics $g_j$ converge to a limit $g_j\to g$ in $C^0$. This
implies that for any fixed metric on $\FFF$ the corresponding
energies of $\phi_j'$ are also uniformly bounded.

At this point we may apply Theorem 1 in \cite{IvSh} on Gromov  
compactness for pseudoholomorphic maps
(which is valid for continuous complex structures
on the target manifold converging uniforming to a limit)  
to deduce the existence (again, after restricting to a subsequence)
of a cusp curve $\Sigma^{\cusp}$
together with a limit map $\phi:\Sigma^{\cusp}\to\FFF$, which is
$I$ holomorphic. Their theorem gives an identification of the principal
component $\Sigma_0\simeq\Sigma$ such that the maps $\phi_j'$ converge
pointwise to $\phi$ in $\Sigma_0\setminus\{x_1,\dots,x_p\}$, where
$\{x_1,\dots,x_p\}$ are the bubbling points of $\Sigma_0$ (i.e., the
points where a bubble gets attached). From this we easily deduce that 
$\phi_0:\Sigma_0\simeq\Sigma\to\FFF$ is a section of the bundle
$\FFF$. Now, using the {\it a priori} estimates given in \S3 of \cite{IvSh},
and repeating the argument in the proof of Theorem \ref{regularitat}
we deduce that for any compact
$K\subset\Sigma_0\setminus\{x_1,\dots,x_p\}$ the restriction of
$(A,\phi_0)$ to $K$ is smooth and we have convergence in
$C^{\infty}$ of $(A_j',\phi_j')$ to $(A,\phi_0)$. This implies that
the restriction of $A$ to $\Sigma_0\setminus\{x_1,\dots,x_p\}$ is in
Coulomb gauge. On the other hand, since $\|A-A_0\|_{L^p_1}<\infty$,
we deduce that $A$ is in Coulomb gauge in the whole $\Sigma$, 
so Theorem \ref{regularitat}
implies that $(A_0,\phi_0)$ is smooth. Finally, statement (3) of the
theorem follows from statement (3) of Theorem 1 in \cite{IvSh}.
\end{pf}

\noindent {\bf Remark.} In \cite{Mu2} a proof of this theorem is given
which applies to the case of compact connected structure group
different from $S^1$ (see the remarks in Section \ref{conditions}). 
The idea consists of first proving a local
result (essentially, combining Uhlenbeck's theorem on existence of
local Coulomb gauge with and {\it equivariant} version of
Gromov--Schwarz lemma, which is Lemma 4.2.1 in \cite{Mu1}), and then
using a standard patching argument.

\section{Invariant complex structures and moduli of rational curves}
\label{s:rationalcurves}

In all this section $\Sigma$ will be the Riemann sphere $S^2=\CP^1$. 
It is well known that for a generic complex structure $I\in\End(TF)$
compatible with $\omega$ the moduli space
of simple $I$ holomorphic maps $s:\Sigma\to F$ is a smooth manifold of
dimension $2\la c_1(TF),s_*[\Sigma]\ra+2n$ 
(see for example Theorem 3.1.2 in
\cite{McDS1}). (Recall that $s$ is a simple map if it does
not factor through a nontrivial ramified covering $\Sigma\to\Sigma$.)

However, in order for the equations (\ref{equs2}) to be gauge
invariant we need to chose $S^1$-invariant complex structres on $F$,
and these are in general far from being generic. In fact, 
most of the times the moduli of simple holomorphic maps with respect
to $S^1$-invariant complex structures will not be smooth or will not
have the expected dimension. So to get smooth moduli of maps we will
have to restrict ourselves to subsets of the set of simple holomorphic
maps. We will take these subsets to be simple curves with fixed
isotropy pair (see below).

The results in this section might be seen as a piece of Gromov--Witten
theory for symplectic orbifolds. Such a theory should study in general
pseudo-holomorphic maps from compact complex orbifolds of complex
dimension 1 to orbifolds, and here we will study in particular
pseudo-holomorphic maps from $\CP^1/\ZZ_m$ to $F/\ZZ_m$, where
$\ZZ_m=\ZZ/m\ZZ$ acts on $\CP^1$ by rotations through a fixed axis.
Parts of Gromov--Witten theory for orbifolds in the algebraic category
have been worked out by D. Abramovich and A. Vistoli \cite{AbVi}.

Let $\III_{\omega}$ be the set of complex structures
on $F$ which are compatible with $\omega$, and let
$\III_{\omega,S^1}\subset\III_{\omega}$ be the $S^1$-invariant ones. 
(Recall that by Lemma 5.49 in \cite{McDS2} $\III_{\omega,S^1}$
is a nonempty and contractible set in the $C^{\infty}$ topology.)

\subsection{Isotropy pairs}
\begin{definition}
Let $s:\Sigma\to F$ be any smooth map. We define the isotropy pair
of $s$ to be the pair of closed subgroups 
$L(s)\subset H(s)\subset S^1$ defined as follows
\begin{align*}
H(s) &:= \{\theta\in S^1|\ \theta\cdot s(\Sigma^1)=s(\Sigma^1)\} \\
L(s) &:= \{\theta\in H(s)|\ \theta|_{s(\Sigma^1)}=\Id\}.
\end{align*}
(The dot $\cdot$ means the action of $S^1$ on $F$.)
\label{defisopairs}
\end{definition}
% Observe that $L(s)$ is the biggest subgroup $\Gamma\subset S^1$
% such that $s(\Sigma)\subset F^{\Gamma}$.

\begin{theorem}
Let $I\in\III_{\omega,S^1}$, and let $s:\Sigma\to F$ 
be a simple holomorphic map. Let $H=H(s)$. There exists
a disk $D\subset\Sigma$ such that
$$S^1\cdot s(D)\cap s(\Sigma)=H\cdot s(D).$$
\label{stabilH}
\end{theorem}
\begin{pf}
In order to prove the theorem we will use the following result
on holomorphic curves (see Lemma 2.2.3 in \cite{McDS1}).
\begin{lemma}
Let $I\in\III_{\omega}$, and let $s_1,s_2:\Sigma\to F$ 
be two simple $I$-holomorphic maps. Let 
$K\subset\Sigma$ be a closed subset such that, for any $x\in K$,
$ds_1(x)\neq 0\neq ds_2(x)$. 
If the intersection $s_1(K)\cap s_2(K)$ contains infinite points, then
$s_1=s_2$.
\label{interseccio}
\end{lemma}
From now on we fix a complex structure $I\in\III_{\omega,S^1}$, and
we take on $F$ the $S^1$-invariant metric $\omega(\cdot,I\cdot)$. 
Let $\fX\in \Gamma(TF)$ be the vector field generated
by the infinitesimal action of
$\imag\in\imag\RR=\Lie(S^1)$. For any $x\in\Sigma$ and
any smooth map $s:\Sigma\to F$ we define
$$\theta_s(x):=\dist(\fX(s(x)),ds(T\Sigma)(x))$$
(note that $\fX(s(x))\in T_{s(x)}F$ and that $ds(T\Sigma)(x)$ is a
2-dimensional subspace of $T_{s(x)}F$).
Suppose to begin with that $H=H(s)=\{1\}$. Assume
that for any open set $U\subset\Sigma$ there exists a point 
$x\in U$ and $1\neq \alpha\in S^1$ such that $\alpha\cdot s(x)\in s(\Sigma)$.
We will see that this leads to a contradiction.

Let $Z=s^{-1}(\{s(z)|z\in\Sigma,\ ds(z)=0\})$. 
This is a finite set (see Lemma 2.2.1 in \cite{McDS1}).
So $s(\Sigma)$ is not contained in $S^1\cdot s(Z)$, since the
latter is a disjoint union of points and circles (since $s$ is simple,
it is in particular non-constant). Let $T$
be a small $S^1$-invariant tubular neighbourhood of $S^1\cdot s(Z)$, and
put $\Sigma'=s^{-1}(F\setminus T)$.
The set of noninjective points $Z'=\{z\in\Sigma|\sharp s^{-1}s(z)>1\}$
can only accumulate at critical points (combine
Lemma 2.2.3 and Proposition 2.3.1 
in \cite{McDS1}), so $Z''=Z'\cap\Sigma'$ is finite.
Hence, $s(\Sigma')$ is not wholly contained in $S^1\cdot s(Z'')$ 
so we may take a small $S^1$-invariant tubular neighbourhood $T''$
of $S^1\cdot s(Z'')$ so that $\Sigma''=s^{-1}(F\setminus T'')\cap\Sigma'$ 
has nonempty interior. 

Let $Y=\{z\in\Sigma|\theta_s(z)=0\}$. This is a closed set.
If $\inter Y\neq\emptyset$ then, for $\alpha\in S^1$ near $1\in S^1$,
$\alpha\cdot s(\Sigma)$ and $s(\Sigma)$ meet at
an open set and hence, by Lemma \ref{interseccio}, they coincide.
But this is implies that $H\neq\{1\}$, in contradiction with our assumption.
So we may suppose that there is a small open disk $D_a\subset\Sigma''$
such that $\inf \theta_s|_{D_a}=a>0$. Suppose also that
$S^1\cdot s(D_a)\subset W\subset F$, where $W$ is open and 
$S^1$-invariant, and all points in $W$ have the same stabiliser,
so that $W/S^1$ is a smooth manifold.

The composition $D_a\stackrel{s_a}{\longrightarrow} W
\stackrel{\pi}{\longrightarrow} W/S^1$
is an embedding (here $s_a=s|_{D_a}$)
if $D_a$ is small enough. Let 
$N\subset W/S^1$ be an open neighbourhood of
$\pi s_a(D_a)$ with a submersion $p:N\to\pi s_a(D_a)$ which is
a left inverse for the inclusion $\pi s_a(D_a)\hookrightarrow N$.
Let $Y_N=Y\cap (\pi s)^{-1}(N)$ (where $\pi:F\to F/S^1$ is the
projection). The critical points of
$$\Sigma\cap (\pi s)^{-1}(N)\stackrel{s}{\longrightarrow}
N\stackrel{p}{\longrightarrow}\pi s(D_a)$$
contain $Y_N$. Hence, by Sard's theorem $\pi s(Y_N)\subset \pi s(D_a)$
has measure zero. Since $\pi s(Y_N)$ is closed, its complementary
contains a closed disk $\Sigma_0$. Furthermore, there exists
$b>0$ such that for any $x\in \Sigma_0$ and $\alpha\in S^1$
if $\alpha\cdot s(x)=s(y)\in s(\Sigma)$, then
$\theta_s(y)\geq b$.

From the construction of $\Sigma_0$ we deduce the following. 
There exist real positive numbers $r,\ \eta,\ \epsilon$ such that
for any $x\in\Sigma_0$ and $\alpha\in S^1$ if
$z=\alpha\cdot s(x)\in s(\Sigma)$, $s^{-1}(z)$ has a unique element
$y\in\Sigma$ and if $D_y=D(y;r)$ is the disk centered at $y$
of radius $r$, the following holds.
\begin{enumerate}
\item[P1.] If $w\in s(\Sigma)$ and $d(w,y)<\eta$, then $w\in s(D_y)$.
\item[P2.] There exists an open neighbourhood $V\subset F$ of
$s(y)$ containing $s(D_y)$ and a chart
$\phi=(\phi_1,\dots,\phi_{2n}):V\to\RR^{2n}$ with 
$\phi(s(y))=0$ such that
\begin{itemize}
\item[P2a.] For any $v\in D_y$, $\phi_3(s(v))=\dots=\phi_{2n}(s(v))$.
\item[P2b.] If $\beta\in[-\epsilon,\epsilon]\subset S^1$, then for
any $v\in D_y$, $\beta\cdot s(v)\in V$ and
$$\phi(\beta\cdot s(v))=\phi(s(v))+(0,0,\beta,0,\dots,0).$$
\end{itemize}
\end{enumerate}

We assume for the rest of the argument that $\diam(s(\Sigma_0))<\eta/2$.
Let us identify $S^1\simeq [0,2\pi)$ so that $0$ is the identity
and consider
\begin{equation}
I=\{(\alpha,x)\in (0,2\pi)\times\Sigma_0|\alpha\cdot s(x)\in s(\Sigma)\}.
\label{defI}
\end{equation}
Thanks to the inequality $\theta_s|_{\Sigma_0}\geq b$ we know that
there exists $\delta>0$ such that
$I\subset [\delta,2\pi-\delta]\times \Sigma_0$. Clearly $I$ is closed.
By our assumption the image of the projection 
$\pi_{\Sigma}:I\to\Sigma_0$ is dense and so (since it is also closed)
coincides with $\Sigma_0$. Let now $[0,\mu]\subset[-\epsilon,\epsilon]$ be 
a subset such that for any $\nu\in[0,\mu]$ and for any 
$x\in F$, $d(x,\nu\cdot x)<\eta/2$.

Cover $[\delta,2\pi-\delta]$ with closed intervals $A_1,\dots,A_r$ of
length $<\mu$ and let $I_k=I\cap A_k\times\Sigma_0$.
Since $\pi_{\Sigma}(I_1)\cup\dots\cup\pi_{\Sigma}(I_r)=\Sigma_0$
and $\pi_{\Sigma}(I_l)$ is closed for any $l$, there exists
a $\pi_{\Sigma}(I_k)$ with nonempty interior. Let $D\subset
\inter \pi_{\Sigma}(I_k)$ be a disk, and take $x\in D$.
By assumption there exists $\alpha\in A_k$ such that
$\alpha\cdot s(x)=s(y)$, $y\in\Sigma$. By P2 there exists
an open set $V\subset F$ containing $s(D_y)=s(D(y,r))$ and a chart
$$\phi:V\to\RR^{2n}.$$
For any $z\in D$ there exists $\beta\in A_k$ such that 
$\beta\cdot s(Z)\in s(\Sigma)$. On the other hand, since
$d(\alpha\cdot s(z),\alpha\cdot s(x))=d(s(z),s(x))<\eta/2$
and $|\alpha-\beta|<\mu$ we have
$$d(\beta\cdot s(z),s(y))<\eta.$$
Hence, by P1, 
$\beta\cdot s(z)\in s(D_y)$. 
So by P2b, if $w=s(z)$, then $\phi_3(w)=\alpha-\beta$,
$\phi_4(w)=\dots=\phi_{2n}(w)=0$. This implies that for any
$z\in D$, $\sharp I_k\cap\{z\}\times A_k=1$.
Let $(z,h(z))$ be the unique element of this set. The function
$h:D\to A_k$ is $h(z)=\alpha-\phi_3(s(z))$ and so is continuous.
Hence there exists $c\in A_k$ such that $\sharp h^{-1}(c)=\infty$
(this follows from this easy result:
if $h:[0,1]^2\to [0,1]$ is a continuous map, then 
there exists $c\in I$ such that $\sharp h^{-1}(c)=\infty$). 
From this we see that 
$c\cdot s(\Sigma)\cap s(\Sigma)$ has infinite points which do not
accumulate on critical points (since $s(\Sigma_0)$ is at positive distance
from the $S^1$-orbit of the image of any critical point of $s$). 
Finally, using Lemma \ref{interseccio} we deduce that
$c\cdot s(\Sigma)=s(\Sigma)$, in contradiction with the
assumption $H=\{1\}$. This finishes the proof of the case $H=\{1\}$.

The case $H=S^1$ is trivial. 
Suppose to finish now that $1<\sharp H<\infty$. We assume that
for any open set $U\subset\Sigma$ there exists $x\in U$ and
$\alpha\in S^1\setminus H$ such that $\alpha\cdot s(x)\in s(\Sigma)$.
We do exactly the
same thing as in the case $H=\{1\}$ to get a subset $\Sigma_0\subset\Sigma$
(note that the function $\theta_s(x)$ is equivariant under the action of
$H$). Now, the set $I$ defined in (\ref{defI}) is at
positive distance from $H\times\Sigma_0$. So the element
$c\in S^1$ found at the end of the reasoning does not belong
to $H$, and hence the fact that $c\cdot s(\Sigma)=s(\Sigma)$ leads to  
a contradiction.
\end{pf}

Let $s:\Sigma\to F$ be a simple map, and let $g\in H(s)$ be any element.
Let $\Sigma_i$ be the set of injective points of $s$, that is,
$\Sigma_i=\{x\in\Sigma|\ ds(x)\neq 0,\ \sharp s^{-1}s(x)=1\}.$
The action of $g$ on $s(\Sigma)$ induces a holomorphic bijection
$\gamma_i(g):\Sigma_i\to\Sigma_i$
which can be extended to a homeomorphism
$\gamma(g):\Sigma\to\Sigma.$
Now, since the map $s$ is simple, 
the noninjective points $\Sigma\setminus\Sigma_i$ can
only accumulate at a finite set of points (namely, the critical
points $\Ker ds$), and hence the map $\gamma(g)$ is holomorphic
by standard removability of singularities. This way we have defined
a map $\gamma:H(s)\to\Aut(\Sigma)=\PSL(2;\CC)$. Obviously, $\Ker\gamma=L(s)$.
Let $$\Map^L(\Sigma,F)=\{s\in\Map(\Sigma,F^L)\mid L(s)=L\}.$$

\begin{theorem}
Let $I\in\III_{\omega,S^1}$,
and let $s\in\Map^L(\Sigma,F^L)$ be a simple holomorphic map.
Then the set $\{x\in\Sigma\mid L\neq (S^1)_{s(x)}\}$  
(where $(S^1)_{s(x)}$ denotes the stabiliser of $s(x)\in F$)
is finite.
\label{stabilL}
\end{theorem}
\begin{pf}
Let $\Sigma'=\{x\in\Sigma\mid L\neq (S^1)_{s(x)}\}$ and
suppose that $\sharp\Sigma'=\infty$.
Since the set of different stabilizers of points of $F$ is finite, we 
may assume that there exists a group $L''$ strictly containing $L$ such that
$$\Sigma''=\{x\in\Sigma\mid (S^1)_{s(x)}=L''\}$$
has infinite elements. Let now $\theta\in L''\setminus L$.
Then $s(\Sigma)$ and $\theta\cdot s(\Sigma)$ intersect 
at an infinite set $\Sigma''$ of points. Hence by Lemma
\ref{interseccio} they coincide, and so $\theta\in H(s)$. 
But now $\gamma(\theta)\in\Aut(\Sigma)$ has infinitely many 
fixed points (all the points in $\Sigma''$), and so it must 
be the identity. But this implies 
that $\theta\in L$, which is a contradiction.
\end{pf}

\subsection{The moduli of rational curves}
Let $L\subset S^1$ be a closed group. 
The fixed point set $F^L\subset F$
is a compact symplectic submanifold (with possibly several
connected components of different dimension). 
The action of $S^1$ on $F$ gives an action of the Lie group
$S^1/L$ on $F^L$.
Fix a closed subgroup $\Gamma\subset\Aut(\Sigma)$, and
assume that there is an injection $\rho:\Gamma\to S^1/L$. 
We will say that a map $s:\Sigma\to F^L$ is
$(\Gamma,\rho)$-equivariant if $s(gx)=\rho(g)s(x)$ for any
$x\in\Sigma$ and $g\in\Gamma$. Let us define
\begin{equation}
\Map(L,\Gamma,\rho)
=\{s\in\Map(\Sigma,F^L)_{L^p_1}\mid L(s)=L,\text{ $s$ is 
$(\Gamma,\rho)$-equivariant}\}.
\end{equation}

Let  $I\in\III_{\omega,S^1}$. Define the moduli of 
$(L,\Gamma,\rho)$-equivariant curves with respect to $I$ to be
$$\MMM_I(L,\Gamma,\rho)
=\{s\in\Map(L,\Gamma,\rho)\mid \ov{\partial}_Is=0,
\text{ $s$ simple}\}.$$
For any $B\in H_2(F;\ZZ)$, let also 
$\MMM_I(L,\Gamma,\rho;B)=
\{s\in\MMM_I(L,\Gamma,\rho)\mid s_*[\Sigma]=B\}.$

\begin{theorem}
There is a subset $\III^{L,\Gamma,\rho}\subset\III_{\omega,S^1}$
of Baire second category (with respect to the $C^{\infty}$ topology on
$\III_{\omega,S^1}$) such that for any $I\in\III^{L,\Gamma,\rho}$
the moduli space $\MMM_I(L,\Gamma,\rho)$ is smooth and oriented. 
Furthermore, for any $I_0,I_1\in\III^{L,\Gamma,\rho}$, there exists a path
$[0,1]\ni\lambda\mapsto I_{\lambda}\in\III_{\omega,S^1}$
such that the space
$$\bigcup_{\lambda\in[0,1]}\MMM_{I_\lambda}(L,\Gamma,\rho)$$
has a natural structure of smooth oriented cobordism between 
$\MMM_{I_0}(L,\Gamma,\rho)$ and $\MMM_{I_1}(L,\Gamma,\rho)$.
\label{modulireg}
\end{theorem}
\begin{pf}
The proof, with due modifications, is exactly like that of 
Theorem \ref{modulipert} or of Theorem 3.1.2 in \cite{McDS1}.
We will give a little more details, since in the course of the proof
one needs to use Theorems \ref{stabilH} and \ref{stabilL}.
We start considering the completion $\III_{\omega,S^1}^l$ (resp.
$\III_{\omega}^l$) of $\III_{\omega,S^1}$ (resp. $\III_{\omega}^l$) 
in the $C^l$ norm, where $l>0$ is a big
enough integer, and we define for any $B\in H_2(F;\ZZ)$
$$\MMM_{\III^l}(L,\Gamma,\rho;B)
=\left\{(s,I)\in\Map(L,\Gamma,\rho)
\times \III_{\omega,S^1}^l\Big|
\begin{array}{l}\ov{\partial}_Is=0,\ s_*[\Sigma]=B,\\
\mbox{and $s$ simple }\end{array}\right\}.$$
We next prove that $\MMM_{\III^l}(L,\Gamma,\rho;B)$
is a smooth Banach manifold.
Let $(s,I)\in\MMM_{\III^l}(L,\Gamma,\rho;B)$
be any point. We have to check that the linearisation
$$D\EEE(u,I):\Omega^0(s^*TF^L)^{\Gamma}_{L^p_1}\times T_I\III_{\omega,S^1}^l
\to \Omega^{0,1}_I(s^*TF^L)^{\Gamma}_{L^p}$$
of the equation at $(s,I)$ is surjective. Here we denote by 
$\Omega^{0,1}_I(s^*TF^L)^{\Gamma}$ the $\Gamma$-invariant sections 
of $\Lambda^{0,1} T\Sigma\otimes_{\CC} s^*TF^L$ (the subscript
stresses the fact that when tensoring over $\CC$ we use the complex
structure $I$; on the other hand, this bundle has an
action of $\Gamma$ through the representation $\rho$).
The tangent space $T_I\III_{\omega,S^1}^l\subset T_I\III_{\omega}^l$ 
is equal to the subspace of $\Gamma$-invariant elements in 
$T_I\III_{\omega}^l$. This latter space is the set of $C^l$ sections
of the bundle $\End(TF,I,\omega)$ whose fibre at $x\in F$ is the
space of linear maps $Y:T_xF\to T_xF$ which satisfy
$$YI+IY=0\mbox{ and }\omega(Y\cdot,\cdot)+\omega(\cdot,Y\cdot)=0$$
(see p. 34 in \cite{McDS1}). 

We now follow the notation (and the ideas) of the proof of Proposition 
3.4.1 in \cite{McDS1}. We may write the differential
$D\EEE(s,I)(\xi,Y)=D_s\xi+\frac{1}{2}Y(s)\circ ds\circ j,$
where $j$ is the complex structure in $\Sigma$ and $D_s$ is a first 
order differential operator whose symbol coincides with that of
Cauchy-Riemann operator. Hence $D_s$ is elliptic and consequently
Fredholm. So if $D\EEE(s,I)$ 
were not exhaustive there would exist a nonzero element
$\eta\in \Omega^{0,1}_I(s^*TF^L)^{\Gamma}_{L^q}$ (where $1/p+1/q=1$) such
that for any $\xi\in \Omega^0(s^*TF^L)^{\Gamma}$ and for any
$Y\in T_I\III_{\omega,S^1}^l$
\begin{equation}
\int_{\Sigma}\la\eta,D_s\xi\ra=0
\mbox{ and }\int_{\Sigma}\la\eta,Y(s)\circ ds\circ j\ra=0.
\label{conseq}
\end{equation}

We now invoque Theorem \ref{stabilH} and obtain a disk $D\subset\Sigma$
such that $$S^1\cdot s(D)\cap s(\Sigma)=H\cdot s(D).$$
Using theorem \ref{stabilL} we deduce that (after possibly
shrinking $D$) all the elements in $s(D)$ have stabiliser equal to $L$.
Then $\eta$ vanishes on an open subset of $D$. For suppose that
$\eta(x)\neq 0$, where $x\in D$. One can always find an endomorphism
$Y_0\in\End(T_{s(x)}F,I_{s(x)},\omega_{s(x)})^L$ such that 
$\la \eta(x),Y_0\circ ds(x)\circ j(x)\ra\neq 0,$
since $\eta(x)\in T_{s(x)}F^L$.
We extend $Y_0$ to $S^1\cdot s(x)$ in a $S^1$-equivariant 
way (we can do this because $(S^1)_{s(x)}=L$ and we took $Y_0$ to 
be $L$-invariant) and then we use a $S^1$-invariant smooth cutoff 
function to extend $Y_0$ to a small neighbourhood of $S^1\circ s(x)$.
This can be done in such a way that the right
hand side integral in (\ref{conseq}) does not vanish. And this is
a contradiction. 

Consequently $\eta$ vanishes in $D$. Since it also satisfies the left 
hand side equation in (\ref{conseq}), Aronszajn's theorem \cite{Ar} 
(see Theorem 2.1.2 in \cite{McDS1}) implies that $\eta$ vanishes identically. 
So $D\EEE(s,I)$ must be exhaustive, and this finishes
the proof that $\MMM_{\III^l}(L,\Gamma,\rho;B)$ is smooth.

The proof of the first statement in 
Theorem \ref{modulireg} is resumed as
Theorem \ref{modulipert} or in p. 36 in \cite{McDS1}. One uses the
Sard--Smale theorem (for that $l$ has to be big enough, depending on the 
index of the linearisation $D\EEE$, which on its turn is a function
of $B\in H_2(F;\ZZ)$) to prove the existence of
a subset $(\III^{L,\Gamma,\rho})^{B,l}\subset\III_{\omega,S^1}^l$ 
of the second category
such that for any $I\in (\III^{L,\Gamma,\rho})^{B,l}$ the moduli space
$\MMM_I(L,\Gamma,\rho;B)$ is smooth. Then a trick of
Taubes allows to deduce from this that
there exists a subset $(\III^{L,\Gamma,\rho})^B\subset\III_{\omega,S^1}$
of the second category with the same property, but consisting
of smooth complex structures and not of $C^l$ ones as before.
Since the set of homology classes $B\in H_2(F;\ZZ)$ is
countable, the intersection 
$$\III^{L,\Gamma,\rho}=\bigcap_{B\in H_2(F;\ZZ)}(\III^{L,\Gamma,\rho})^B$$
is again of the second category.

To finish the proof, note that the linearisation of the equations
is, modulo a compact operator, the Cauchy--Riemann operator. Hence
the cohomology groups of the deformation complex carry natural
orientations (because they are complex vector spaces) and 
consequently so does the moduli space.

The last statement of the theorem on the independence of the
cobordism class for generic $I$ is proved analogously.
\end{pf}

We now define 
$$\III^{\reg}_{\omega,S^1}=\bigcap_{L,\Gamma,\rho} \III^{L,\Gamma,\rho},$$
where the intersection is taken for the triples $(L,\Gamma,\rho)$
such that the moduli $\MMM_{\III}(L,\Gamma,\rho;B)$ is nonempty
for some $B\in H_2(F;\ZZ)$. Again, this is a Baire set of the
second category.

\subsection{Index computations}
In this subsection we will compute the dimension of the moduli
spaces $\MMM_I(L,\Gamma,\rho)$ for generic $I$.
Let us fix a triple $(L,\Gamma,\rho)$. Recall that $L$ is a subset of $S^1$,
$\Gamma$ is a compact subgroup of $\Aut(\CP^1)$ and $\rho:\Gamma\to S^1/L$
is an injection. This latter condition implies that $\Gamma$ is
abelian. Since $\Aut(\CP^1)=\PSL(2;\CC)$, any element in 
$\Aut(\CP^1)$ which spans a compact subgroup must fix two points
of $\CP^1$. And since $\Gamma$ is abelian, there
must exist two points $x_+$ and $x_-$ which are fixed by all the
elements of $\Gamma$. Using one of the fixed points, say $x_+$, we
get an injection $\Gamma\to S^1\subset\CC^*$ by assigning to
any $\gamma\in\Gamma$ the induced endomorphism 
$\iota(\gamma)\in\GL(T_{x_+}\CP^1)$.
In the sequel we will identify $\Gamma$ with its image in $S^1$.
There are two possibilities. Either $\Gamma$
is a finite group or $\Gamma\simeq S^1$. When $\Gamma$ is a finite group,
the map $\iota$ fixes an isomorphism $\Gamma\simeq\ZZ/m\ZZ$, and
when $\Gamma$ is infinite $\iota$ gives an identification with $S^1$.

If $\Gamma\neq\{1\}$ then, for any $(\Gamma,\rho)$-equivariant
map $s:\CP^1\to F$, the fixed points $x_\pm$ are mapped by $s$
to the fixed point set $F^f$ (because by assumption
the action on $F\setminus F^f$
is free). Let $z$ be a holomorphic coordinate in $\CP^1$ centered at
$x_+$. Taking $S^1$-equivariant coordinates in a neighbourhood of $s(x_+)$ 
the map $s$ can be written (see p. 16 in \cite{McDS1})
$s(z)=az^l+O(|z|^{l+1})$,
and the constant $a$ can be identified with an element of $T_{s(x_+)}F$.
Let $T_{x_\pm}^PF$ (resp. $T_{x_\pm}^ZF$, $T_{x_\pm}^NZ$) be the subspace
of $T_{x_\pm}F$ spanned by vectors of weight $1$ (resp. $0$, $-1$)
under the action of $S^1$. Since the action of $S^1$ on $F$ is
semi-free, there are no more weights, and hence $a$ must
lie in $T_{x_\pm}^PF\cup T_{x_\pm}^ZF\cup T_{x_\pm}^NF$ (otherwise
the vector space it spans 
would not be invariant under the action of $\Gamma$). Using
the local expression of $s(z)$ 
we may write, for any $\theta\in\Gamma$ and $z$ near $x_+$,
$s(\theta z)=\rho(\theta)\cdot s(z)=\theta^l s(z)$
modulo $O(|z|^{l+1})$. The $\cdot$ in the second term refers to
the action of $S^1$ on $T_{x_+}F$. From this we deduce that
$a$ cannot belong to $T_{x_\pm}^ZF$ and that:
if $a\in T_{x_\pm}^PF$ then $\rho(\theta)=\theta^l$, and 
if $a\in T_{x_\pm}^NF$ then $\rho(\theta)=\theta^{-l}$.
In fact, after possibly composing $s$ with the holomorphic
map $r:\CP^1\to\CP^1$ defined $r([x:y])=[y:x]$ in coordinates
for which $x_+=[0:1]$ and $x_-=[1:0]$, we may assume
that $a\in T_{x_\pm}^PF$. Hence $\rho(\theta)=\theta^l$ for
any $\theta\in\Gamma$, where $l$ is a positive integer.
If $\Gamma=S^1$, then $l$ must be $1$, and if $\Gamma=\ZZ/m\ZZ$
then $l$ and $m$ must be coprime and the representation $\rho$
only depends on the class of $l$ modulo $m$. 

In the sequel we will write $\MMM_I(L,\Gamma,l;B)$ instead of
$\MMM_I(L,\Gamma,\rho;B)$. When $\Gamma=1$ we will write 
$\MMM_I(L;B)$ instead of $\MMM_I(L,\Gamma,\rho)$,
and when $L=\Gamma=1$ we will write $\MMM_I(B)$.

\subsubsection{The deformation complex}
\label{defocplx}
Let $B\in H_2(F;\ZZ)$, $I\in\III^{\reg}_{\omega,S^1}$, and
$s\in\MMM_I(L,\Gamma,l;B)$. The deformation complex
of the moduli $\MMM_I(L,\Gamma,l;B)$ at $s$ is
$$D_s^{\Gamma}:\Omega^0(s^*TF^L)^{\Gamma}\to
\Omega^{0,1}_I(s^*TF^L)^{\Gamma},$$
where $D_s^{\Gamma}$ is equal to the Cauchy-Riemann operator modulo a compact
operator (see p. 28 in \cite{McDS1}). Since $I\in\III^{\reg}_{\omega,S^1}$,
this operator is exhaustive and consequently the dimension of 
$\MMM_I^{L,\Gamma,l}(B)$ at $s$ is equal to 
$\dim (\Ker D_s^{\Gamma})$. To compute this dimension we consider the 
natural extension of $D_s^{\Gamma}$
\begin{equation}
D_s:\Omega^0(s^*TF^L)\to\Omega^{0,1}_I(s^*TF^L)
\label{defcomg}
\end{equation}
(this is the deformation complex of the moduli of holomorphic
curves in $F^L$). The operator $D_s$ is $\Gamma$-equivariant,
and hence acts on the cohomology groups $H^i_s$ of the complex.
We have $\Ker D_s^{\Gamma}=(H^0_s)^{\Gamma}$
and $\Coker D_s^{\Gamma}=(H^1_s)^{\Gamma}=0$. So the complex
dimension at $s$ is equal to 
\begin{equation}
\dim T_s\MMM_I(L,\Gamma,l;B)=\dim (H^0_s)^{\Gamma}
-\dim (H^1_s)^{\Gamma}.
\label{ladimensio}
\end{equation}
This dimension can be computed putting instead of $D_u$ any
equivariant Dolbeaut operator 
on $s^*TF^L$, since they have the same symbol. Because the action of $S^1$ 
on $F$ is almost-free, we need only distinguish these possibilities.

\noindent{\bf Case 1.} $L=S^1$, $\Gamma=\{1\}$. 
Let $F^f=F_1\cup\dots\cup F_r$ be the connected 
components of the fixed point set. Suppose that 
$B\in H_2(F_k;\ZZ)\subset H_2(F;\ZZ)$. Then by Riemann-Roch
the moduli space has dimension
$$\dim \MMM_I(S^1;B)=2\la c_1(TF_k),B\ra+\dim F_k.$$

\noindent{\bf Case 2.} $L=\{1\}$, $\Gamma\neq\{1\}$. 
Let $x_{\pm}\in\CP^1$ be the two points which are fixed by $\Gamma$.
Since the map $s$ is $(\Gamma,\rho)$-equivariant, we have a natural lift of 
the action $\rho$ of $\Gamma$ on $\CP^1$ to $E=s^*TF\to\CP^1$. 
Let us write it $\gamma:\Gamma\to\Aut(E)$,
where $\Aut(E)$ denotes the automorphisms of $E$ as vector bundle.
The map $\gamma$ induces representations $\gamma_{\pm}$ of 
$\Gamma$ on the fibres $E_{x_\pm}$ over $x_{\pm}$. The weights of this 
representation are $l$ times the weights of the representation of $S^1$ 
on $TF_{s(x_\pm)}$ (which belong to $\{-1,0,1\}$). 
Let $P_\pm$ (resp. $Z_\pm$, $N_\pm$) be the number of weights
of the representation $\gamma_\pm$ which are equal
to $1$ (resp. $0$, $-1$). 

Denote $\Ind_\gamma(E)=\dim (H^0_s)^{\Gamma}-\dim (H^1_s)^{\Gamma}$ 
the $\Gamma$-invariant part of the index of the operator $D_s$ on $E$.  
We will denote by $\rk(E)$ the complex rank of $E$.

\subsubsection{Case $\Gamma=\ZZ/m\ZZ$}
\begin{theorem}
Let $P_\pm$ (resp. $Z_\pm$, $N_\pm$) be the number of weights
of $\gamma_\pm$ which are $l$ (resp. $0$, $-l$). Let $l'=l+km$ for $k\in\ZZ$
such that $1\leq l'\leq m-1$. 
Then $$\Ind_\gamma(E)=\frac{2}{m}(\deg(E)+m\rk(E)-m(P_-+N_+)+
l'(P_-+N_+-P_+-N_-)).$$
\label{indtor}
\end{theorem}
\begin{pf}
We may write
\begin{align*}
m\Ind_{\gamma}(E) &=
\sum_{k\in\ZZ/m\ZZ}\Tr(\gamma(k),H^0(E))-\Tr(\gamma(k),H^1(E)) \\
&=2(\deg(E)+\rk(E))+
\sum_{k=1}^{m-1}\Tr(\gamma(k),H^0(E))-\Tr(\gamma(k),H^1(E)), 
\end{align*}
by Riemann--Roch,
where $\Tr(\gamma(k),H^i(E))$ denotes the trace of $\gamma(k)$ 
acting on $H^i(E)$. We will compute the value of
$$\Tr(\gamma(k),H^0(E))-\Tr(\gamma(k),H^1(E))$$ 
for $1\leq k\leq m-1$ using Atiyah-Bott fixed point theorem 
(see \cite{BeGeV}).

\begin{theorem}[Atiyah-Bott]
Let $M$ be a compact complex manifold and $W\to M$ a holomorphic
vector bundle. Let $g:M\to M$ be a complex diffeomorphism
which lifts to $g:W\to W$. Suppose that the fixed points of
$g$ are isolated. Then
$$\sum_i (-1)^i\Tr_{\CC}(g,H^i(W))=
\sum_{x_0\in M^g} \frac{\Tr_{\CC}(g_{x_0}^W)}
{\det_{T_{x_0}^{1,0}M}(1-g_{x_0}^{-1})},$$
where $H^i(W)$ is the $i$-th Dolbeaut cohomology group and 
$g_{x_0}^W:W_{x_0}\to W_{x_0}$ is the complex linear endomorphism of
the fibres over the fixed points induced by $g$ (we use the
determinant of complex endomorphisms).
\label{AtiBot}
\end{theorem}

In our case we have for any $1\leq k\leq m-1$ a complex diffeomorphism
$\rho_m(k)\in\Aut(\CP^1)$ whose fixed points are $x_\pm$. Let
$\theta=\exp(2\pi\imag/m)$. We then have
$$\det(1-\rho_m(k)_{x_\pm}^{-1})=(1-\theta^{\mp 1}).$$
Let $N=\rk(E)$ and let $b^1_\pm,\dots,b^N_\pm\in\ZZ/m\ZZ$ be
the weights of $\gamma_\pm$. Then
$$\Tr_{\CC}(\gamma(k)_\pm)=\sum_{j=1}^N\theta^{b^j_{\pm}k}.$$
So using Theorem \ref{AtiBot} we conclude that
\begin{equation}
\Ind_{\gamma}(E)=\frac{2}{m}
\left(\deg(E)+\rk(E)+\sum_{k=1}^{m-1}\sum_{j=1}^N
\left(\frac{\theta^{b^j_+k}}{1-\theta^{-k}}
+\frac{\theta^{b^j_-k}}{1-\theta^k}\right)\right).
\label{quasifi}
\end{equation}
\begin{lemma} Let $\theta=\exp(2\pi\imag/m)$. Then for $1\leq w\leq m-1$
\begin{align*}
\sum_{k=1}^{m-1}\frac{1}{1-\theta^k}&=
\sum_{k=1}^{m-1}\frac{1}{1-\theta^{-k}}=
\frac{m-1}{2} \\
\sum_{k=1}^{m-1}\frac{\theta^{wk}}{1-\theta^k}&=
\sum_{k=1}^{m-1}\frac{\theta^{-wk}}{1-\theta^{-k}}=
-\frac{m-1}{2}+w-1.
\end{align*}
\end{lemma}
\begin{pf}
Let $f(x)=\prod_{k=1}^{m-1}(x-\theta^k)$. We have
$f(x)=1+x+\dots+x^{m-1}$ and
$$\sum_{k=1}^{m-1}\frac{1}{1-\theta^k}=
\frac{f'(1)}{f(1)}=\frac{m(m-1)/2}{m}=\frac{m-1}{2}.$$
In general, for any $1\leq w\leq m-1$ 
\begin{align*}
\sum_{k=1}^{m-1}\frac{\theta^{wk}}{1-\theta^k} &=
\sum_{k=1}^{m-1}\left(-\frac{1-\theta^{wk}}{1-\theta^k}
+\frac{1}{1-\theta^k}\right)
=\sum_{k=1}^{m-1}-(1+\theta^k+\dots+\theta^{(w-1)k})+\frac{m-1}{2}\\
&=-\frac{m-1}{2}+w-1,
\end{align*}
since, for any $w\in \ZZ$, $\sum_{k=1}^{m-1}\theta^{wk}$
is $m-1$ if $m\mid w$ and $-1$ otherwise.
\end{pf}

Now, combining the above lemma with (\ref{quasifi}) we get 
$$\Ind_\gamma(E)=\frac{2}{m}(\deg(E)+m\rk(E)-m(P_-+N_+)+
l'(P_-+N_+-P_+-N_-)),$$
which is what we wanted to prove.
\end{pf}

\subsubsection{Case $\Gamma=S^1$}
\begin{theorem}
Let $P_\pm$ (resp. $Z_\pm$, $N_\pm$) be the number of weights
of $\gamma_\pm$ which are $1$ (resp. $0$, $-1$).
Then $$\Ind_\gamma(E)=2(\rk(E)-(P_-+N_+)).$$
\label{inds}
\end{theorem}
\begin{pf}
For any $p\in\NN$, let $\Gamma_p=\ZZ/2^p\ZZ\subset S^1$, and 
consider the action $\gamma_p:\Gamma_p\to\Aut(E)$ induced by
$\gamma$. Then we clearly have
$\Ind_\gamma(E)=\lim_{p\to\infty}\Ind_{\gamma_p}(E)$, and the 
equality follows then from Theorem \ref{indtor}.
\end{pf}

Using the formula $\deg(E)=P_++N_--P_--N_+$, the above index can also
be written $\Ind_\gamma(E)=2(\deg(E)+\rk(E)-(P_++N_-)).$

\subsubsection{An inequality}

\begin{lemma}
Assume that $\Gamma\neq\{1\}$ and that $\deg(E)> 0$. Then we have
$$\Ind_\gamma(E)\leq 2(\deg(E)+\rk(E))-4.$$
\label{destor}
\end{lemma}
\begin{pf}
Suppose to begin that $\Gamma=\ZZ/m\ZZ$ and that
$1\leq l'\leq m-1$ is as in Theorem \ref{indtor}.  
Since $\Ind_{\gamma}(E)$ is an even integer (because both $H^0_s(E)$ and
$H^1_s(E)$ are complex spaces and the action of $\Gamma$ respects the
complex structure, hence $H^0_s(E)^\Gamma$ and $H^1_s(E)^\Gamma$ are
both complex spaces), it is enough for our purposes to prove
$$\frac{1}{2}\Ind_{\gamma}(E)\leq \deg(E)+\rk(E)-(1+1/m).$$
Writing the value of $\Ind_{\gamma}(E)$ given by Theorem \ref{indtor},
multiplying by $m$ and simplifying we arrive at the (equivalent) inequality
$$m+1\leq (m-1)\deg(E)+(m-l')(P_-+N_+)+l'(P_++N_-),$$
which is a consequence of $P_++N_++P_-+N_-\geq 2$ and $\deg(E)\geq 1$, taking
into account that $m-l'\geq 1$ and $l'\geq 1$.

The case $\Gamma=S^1$ can be deduced from the previous one using the
same limit trick as in the proof of Theorem \ref{inds}.
\end{pf}

\subsection{Evaluation maps are submersions}
\label{submersions}
Let $B\in H_2(F;\ZZ)$.
In this subsection we will generalise the result in \S 6.1
of \cite{McDS1} for curves in $\MMM_I(L,\Gamma,\rho;B)$. 
For any $x\in\CP^1$ we have an evaluation map
$$\ev_x:\MMM_{\III}=\MMM_{\III}(L,\Gamma,\rho;B)\to F$$
which sends any $s\in\MMM_{\III}$ to $\ev_x(s)=s(x)$. When
$\Gamma=1$ theorem 6.1.1 in \cite{McDS1} says that the map $\ev_x$
is a submersion. When $\Gamma\neq 1$ this need not hold any longer.
In fact, we must distinguish two possibilities. If $x\neq x_\pm$,
then the map $\ev_x:\MMM_{\III}\to F$ is a submersion, and if
$x=x_\pm$ then the evaluation map $\ev_x$ takes values in $F^f$
and the map $\ev_x:\MMM_{\III}\to F^f$ is a submersion.
We state this in the following lemma.

\begin{lemma}
Suppose that $\Gamma\neq S^1$. Given $I\in\III_{\omega,S^1}$, a curve
$s\in\MMM_I(L,\Gamma,\rho;B)$ and a point $x\in\CP^1$
different from $x_\pm$ (resp. equal to $x_\pm$) 
there exists $\delta>0$ such that for any $v\in T_{s(x)}F^L$
(resp. for any $v\in T_{s(x)}F^f$) and every 
$0<\rho<r<\delta$ there exists a smooth $\Gamma$-equivariant
vector field $\xi\in\Omega^0(s^*TF^L)^{\Gamma}$ and an infinitesimal
variation of almost complex structure $Y\in T_I\III_{\omega,S^1}$
(see Theorem \ref{modulireg}) such that the following
holds

i) $D_s\xi+\frac{1}{2}Y(s)\circ ds\circ j=0$ (that is, the pair
$(\xi,Y)$ belongs to $T_{(s,I)}\MMM_{\III}$),

ii) $\xi(x)=v$ and

iii) $\xi$ is supported in $\Gamma\cdot B_\delta(x)$ and $Y$ is supported
in and arbitrarily small neighbourhood of 
$s(\Gamma\cdot(B_r(x)\setminus B_{\rho}(x))$.
\label{submersioeq}
\end{lemma}
\begin{pf}
Since the proof is almost the same as that of Lemma 6.1.2 in
\cite{McDS1}, we will just give a sketch and mention the differences.
The first thing to do is to find a local solution $\xi_0$ of
$D_s\xi=0$ in $B_{\delta}(x)$ satisfying $\xi_0(x)=v$. 
This is done by solving a boundary value problem (see proposition
4.1 in \cite{McD1} and the references therein). Then one multiplies
$\xi$ by a cutoff function with support in a neighbourhood of
$B_r(x)\setminus B_{\rho}(x)$ to extend $\xi_0$ to a section
of $s^*TF^L$. One then averages $\xi_0$ by the action of $\Gamma$
and obtains a section $\xi\in\Omega^0(s^*TF^L)^\Gamma$.
Finally, one must modify $I$ by a suitable infinitesimal
$Y\in T_I\III_{\omega,S^1}$ so that {\it i)} is satisfied
(in order to take $Y$ $\Gamma$-equivariant one needs to be careful with
the fixed point locus of the action of $S^1$; this may be done 
using theorem \ref{stabilL}, as was done in the proof
of Theorem \ref{modulireg}). This $Y$ can be taken
fulfilling property {\it iii)}, repeating the argument in 
\cite{McDS1} but taking into account $\Gamma$-equivariance.
\end{pf}

\begin{definition}
We will say that a point $x\in\CP^1$ is critical with respect to the
tuple $(L,\Gamma,\rho)$ if either $L=S^1$ or $\Gamma\neq 1$
and $x=x_\pm$.
\label{defcritic}
\end{definition}

\section{Definition of the invariants}
\label{s:definition}

\subsection{Pseudo-cycles}
In this subsection we will review some basic facts about
pseudo-cycles. All the results which we will use are taken from
\S 7.1 in \cite{McDS1}. 

Let $X$ be a smooth compact $m$-dimensional manifold, and let
$R\subset X$ be a subset. We will say that $R$ has dimension at most
$k$ (and write $\dim R\leq k$) 
if $R$ is contained in the image of a smooth map $g:W\to X$, where $W$
is a $\sigma$-compact $k$ dimensional smooth manifold (recall that a
space is $\sigma$-compact if it can be covered by countably many
compact sets). 

Given a map $f:M\to X$ we define the boundary $\Omega_f$ of $f(M)$ to be
$$\Omega_f=\bigcap_{K\subset M}\ov{f(M\setminus K)},$$
where the intersection runs over all the compact subsets $K$ of
$M$. The set $\Omega_f\subset X$ coincides with the set of all points
in $X$ which are limit of sequences $f(m_j)$, where $m_j$ has no
convergent subsequence in $M$.

\begin{definition}
A $k$-dimensional pseudo-cycle is a smooth map $f:M\to X$, where $\dim
M=k$ and $M$ is oriented, 
such that $\Omega_f\leq k-2$. Two $k$-dimensional pseudo-cycles
$f_0:M_0\to X$ and $f_1:M_1\to X$ are bordant if there exists an
oriented cobordism $W$ between $M_0$ and $M_1$ and a smooth map
$F:W\to X$ extending $f_0$ and $f_1$ such that $\dim\Omega_F\leq k-1$.
\end{definition}

Given two pseudo-cycles in $X$ of complementary dimension there is a well
defined intersection number between them (which coincides with the
usual one if they are cycles, see Lemma 7.1.3 in \cite{McDS1}). 
On the other hand, any homology class in $X$ can be represented by a
pseudo-cycle (see Remark 7.1.1 in \cite{McDS1}). Using these two
facts, one can define a pairing between pseudo-cycles and homology
classes in $X$. The following is Lemma 7.1.4 in \cite{McDS1}.

\begin{lemma}
Every $k$-dimensional pseudo-cycle $f:M\to X$ defines canonically a
map $$\Psi_f:H_{m-k}(X;\ZZ)\to\ZZ.$$
Furthermore, if $f$ and $f'$ are bordant, then $\Psi_f=\Psi_{f'}$.
\label{aplihom}
\end{lemma}

\subsection{Construction of a finite dimensional target}
\label{findimtarget}
The strategy in \cite{McDS1} to define Gromov--Witten invariants
consists of proving that the evaluation map (which takes values in the
compact symplect manifold) is a pseudo-cycle. Our aim is to follow the same
idea. 

Suppose that $\beta\in H_2^{S^1}(F;\ZZ)$,
$c\in\imag\RR\setminus\CCC$, and $\rho\in\PPP_c^{\reg}$ have been fixed.
Let $\bM=\bM_{\rho}(\beta,c)$ and $\MMM=\MMM_{\rho}(\beta,c)$.
Take a positive integer $p$. To give a rigorous sense to formula
(\ref{definv}) when $\MMM$ is not compact, we proceed as follows.
Let 
$$\bEE^p=\AAA\times\SSS\times E^p\to\AAA\times\SSS\times\Sigma^p$$ 
be the $(S^1)^p$-bundle pullback of $E^p\to\Sigma^p$ by the projection
$\pi_{\Sigma}:\AAA\times\SSS\times\Sigma^p\to\Sigma^p$.
This bundle has a natural action of $\GGG$, and its restriction to
$\bM\times\Sigma^p$ descends to give a $(S^1)^p$-bundle
$\EE^p\to\MMM\times\Sigma^p$. Furthermore, the canonical section
$$\begin{array}{rcl}
\obPhi:\AAA\times\SSS\times\Sigma^p & \to & 
\bEE^p\times_{(S^1)^p}F^p=\AAA\times\SSS\times F^p \\
(A,\phi,x_1,\dots,x_p) & \mapsto &
(A,\phi,\phi(x_1),\dots,\phi(x_p))
\end{array}$$ 
is $\GGG$-equivariant and hence
descends to give a section $\bPhi$ of the bundle
$\EE^p\times_{(S^1)^p}F^p$. Now, Lemma \ref{bijeccio} allows to obtain
from $\bPhi$ a map
$$\ev_{S^1}^p:=\bPhi_{\EE^p}:
\MMM\times\Sigma^p\to(F^p)_{(S^1)^p}=(F_{S^1})^p.$$
We would like to treat $\ev_{S^1}^p$ as a pseudo-cycle. 

However, $(F_{S^1})^p$ is not a finite dimensional manifold.
To solve this problem, we will construct a 
compact oriented smooth manifold $T$ with a
$(S^1)^p$ principal bundle $E_T^p\to T$ and a Cartesian diagram
\begin{equation}
\xymatrix{\EE^p\ar[r]^{s^T}\ar[d] & E_T^p\ar[d] \\
\MMM\times\Sigma^p \ar[r] & T.}
\label{modelfinit}
\end{equation}
We will call such a diagram a smooth compact model of $\EE^p$.
Let $$s^T_{F^p}:\EE^p\times_{(S^1)^p}F^p\to E_T^p\times_{(S^1)^p}F^p$$
be the induced map. We will prove that
$s^T_{F^p}\Phi:\MMM\times\Sigma^p\to E_T^p\times_{(S^1)^p}F^p$ is a
pseudo-cycle. 
That this is consistent with our definition, and that the result is
independent of the approximation $E_T^p\to T$ is proved by
formula (\ref{holabondia}) in the remark after Lemma
\ref{bijeccio}.

In the rest of this subsection we will construct the smooth compact
model 
$E_T^p\to T$ and in the next one we will prove that $s^T_{F^p}\Phi$ is
a pseudo-cycle.

\subsubsectionr{}
For any finite subset $P=\{p_1,\dots,p_N\}\subset\Sigma$, 
we let $F_P=\prod_{p\in P}\FFF_p$, 
and we define the map
$$\begin{array}{rcl}
e_P:\AAA\times\SSS & \to & F_P \\
(A,\phi) & \mapsto & (\phi(p_1),\dots,\phi(p_N)).
\end{array}$$ 
Consider the action of
$S^1$ on $F_P$ induced by the action of the constant gauge
transformations on $\FFF$, and let $F_P^f$ denote the fidxed point set
of this $S^1$-action (in fact $F_P\subset\SSS$ supports an action
of the full gauge group $\GGG$). We have $F_P^f=\prod_{p\in P}
\FFF_p^f$, where $\FFF^f=E\times_{S^1}F^f$.

\begin{lemma}
One can take $P\subset\Sigma$ such that 
$$\ov{e_P(\bM)}\cap F_P^f=\emptyset.$$
\end{lemma}
\begin{pf}
For any $\epsilon>0$ we will denote $P_\epsilon\subset\Sigma$ 
any finite subset such that the disks of radius $\epsilon$
centered at the points $p\in P_\epsilon$ cover $\Sigma$. 
Suppose that the claim of the lemma is not true. 
Then there exists a sequence $\epsilon_j\to 0$, sets $P_{\epsilon_j}$
and $\rho$-THCs $(A_j,\phi_j)\in\bM_{\rho}(\beta,c)$ so that for any $j$
the image of the points in $P_{\epsilon_j}$ by the section $\phi_j$
is contained in $\FFF^f=E\times_{S^1}F^f$. By the compactness
Theorem \ref{compactificacio} one may take a subsequence of $(A_j,\phi_j)$
which, after suitably regauging, converge pointwise to a cusp
$\rho$-THC. Now, by construction, the image of the principal component
$\Sigma_0$ of this limit cusp must be inside $\FFF^f$. But this
is in contradiction with our assumption that $\rho\in\PPP_c$.
\end{pf}

Let us take a subset $P\subset\Sigma$ as given by the preceeding lemma.
Let $N$ be a $\GGG$-invariant
tubular neighbourhood of the fixed point set $F_P^f$
which does not meet the closure of $e_P(\bM_{\rho}(\beta,c))$.
Take $S$ to be two
copies of $F_P\setminus N$ glued along $\partial N$:
$$S=(F_P\setminus N)\cup_{\partial N}-(F_P\setminus N).$$
Then $S$ supports an action of $\GGG$, and no point in $S$ is
fixed by a nontrivial constant gauge transformation. Furthermore,
the map $e_P$ gives a $\GGG$-equivariant
map $s_P:\bM_{\rho}(\beta,c)\to S$.

\subsubsectionr{}
Let $d$ be the degree of $E$, and let 
$$\AAA_{\Fl}=\{A\in\AAA\mid F_A=-\imag 2d\pi
\omega_{\Sigma}/\Vol(\Sigma)\}$$ be
the set of projectively flat connections (here $\omega_{\Sigma}$ is 
the symplectic form in $\Sigma$). Then
$$\Jac_d(X):=\AAA_{\Fl}/\GGG$$ is a torus of real dimension 
twice the genus of $\Sigma$. We will construct a retraction 
$\AAA/\GGG\to\Jac_d(\Sigma)$. Recall that we have a metric on $\Sigma$, 
which induces metrics on the exterior 
algebra of forms $\Omega^*(\Sigma)$. Let $\HHH^j$ be the space of harmonic 
$j$-forms with respect to this metric.

Let $$G:\AAA/\GGG\to\Omega^2(\imag\RR)$$
be the map which sends any $[A]$ to
$F_A+\imag 2\pi d\omega_{\Sigma}/\Vol(\Sigma)$.
It is easy to see, using Hodge theory, that the image of $G$ is the 
orthogonal of $\imag\HHH^2(\Sigma)$ in $\Omega^2(\imag\RR)$.
The preimage of $0\in\Omega^2(\imag\RR)$ is precisely $\Jac_d(\Sigma)$. 
In fact, $G:\AAA/\GGG\to\imag\HHH^2(\Sigma)^{\bot}$ is a smooth fibration
with fibres diffeomorphic to $\Jac_d(\Sigma)$. We will construct a connection
on this fibration by specifying its horizontal distribuition.

Given any $[A]\in\AAA/\GGG$, the tangent space $T_{[A]}\AAA/\GGG$
can be canonically identified with $\Ker d_1^*$, where
$d_1:\Omega^0(\imag\RR)\to\Omega^1(\imag\RR)$ is the exterior derivation.
Then we set the horizontal space at $[A]$ to be
$$(T_{[A]}\AAA/\GGG)_h:=\Ker d_1^*\cap (\Ker d_2)^{\bot},$$
where $d_2:\Omega^1(\imag\RR)\to\Omega^2(\imag\RR)$
is the exterior derivation.
Now, using parallel transport along lines going through 
$-\imag 2\pi d\omega_{\Sigma}/\Vol(\Sigma)\in\Omega^2(\imag\RR)$ we get
the desired retraction
$$R:\AAA/\GGG\to\Jac_d(\Sigma).$$

\subsubsectionr{}
By the definition of $S$, the group $\GGG$ acts freely on $\AAA\times
S$ (since the stabiliser of any connection is the set of constant
gauge transformations). We define 
$$T:=(\AAA_{\Fl}\times S)/\GGG\times\Sigma^p.$$ 
To define the bundle $E_T^p\to T$ we do the
following construction.
Let 
$$\bE_S^p=\AAA\times S\times E^p\to\AAA\times S\times\Sigma^p,$$ 
and let $E_S^p=\bE_1^p/\GGG\to(\AAA\times S)/\GGG\times\Sigma^p$ be the
quotient bundle (that $\bE_S^p$ descends follows from the fact that
the action of $\GGG$ on $\AAA\times S$ is free). Let
$j_1:\Jac_d(\Sigma)\to\AAA/\GGG$ and
$j_2:T\to(\AAA\times S)/\GGG\times\Sigma^p$ be the inclusions, and
set $E_T^p:=j_2^*E_S^p$. Let $j_3:E_T^p\to E_S^p$ be the inclusion.
We have the following diagram:
$$\xymatrix
{\bE^p_S \ar[r]\ar[d] & E^p_S \ar[d] & E^p_T\ar[l]_{j_3}\ar[d] \\
\AAA\times S\times\Sigma^p\ar[r]\ar[d] &
(\AAA\times S)/\GGG\times\Sigma^p\ar[d] &
T\ar[d]\ar[l]_-{j_2}\\
\AAA\ar[r] & \AAA/\GGG & \Jac_d(\Sigma)\ar[l]_{j_1}.}$$
The arrows in the second column are fibrations. Taking a connection on
the total fibration $E^p_S\to\AAA/\GGG$ we extend the retraction
$R$ to obtain retractions $R_2$ and $R_3$
of the inclusions $j_2$ and $j_3$. Then we
get a Cartesian diagram
$$\xymatrix{E^p_1\ar[r]^{R_3}\ar[d] & E^p_T\ar[d]\\
(\AAA\times S)/\GGG\times\Sigma^p\ar[r]^-{R_2} & T.}$$
On the other hand, we have a $\GGG$-equivariant map
$$\begin{array}{rcl}
\bM\times\Sigma^p & \to & \AAA\times S\times\Sigma^p \\
(A,\phi,x_1,\dots,x_p) & \mapsto &
(A,\phi(p_1),\dots,\phi(p_N),x_1,\dots,x_p),
\end{array}$$ 
and similarly a lift
$\bM\times E^p\to\AAA\times S\times E^p$.
Dividing out by the action of $\GGG$ and composing with $R_3$ and
$R_2$ we get the desired Cartesian diagram (\ref{modelfinit}).

\subsection{$s^T_F\circ\Phi$ is a pseudo-cycle}

There are now two things to prove. The first one is that, for a
generic choice of complex structure and perturbation $\rho$, 
the map which we obtain using the
finite dimensional approximation (\ref{modelfinit}) is a
pseudo-cycle. The second one is that the bordism class of this
pseudo-cycle is independent of the complex structure and the
perturbation, and that it only depends on the connected component of
$\imag\RR\setminus\CCC$ in which $c$ lies.

\begin{theorem}
(i) Let $\beta\in H_2^{S^1}(F;\ZZ)$ and $c\in\imag\RR\setminus\CCC$.
Let $\QQQ=\PPP_c\times\III_{\omega,S^1}$. 
There is a subset $\QQQ^{\reg}\subset\QQQ$ of Baire of
the second category (with respect to the $C^{\infty}$ topology on
$\QQQ$) such that if $(\rho,I)\in\QQQ^{\reg}$,
$\MMM=\MMM_{\rho,I}(\beta,c)$,
and $s^T_F:\EE^p\times_{(S^1)^p}F^p\to
E^p_T\times_{(S^1)^p}F^p$ is the map induced by $s^T$, then
$s^T_F\circ\Phi:\MMM\times\Sigma^p\to E^p_T\times_{(S^1)^p}F^p$ is a
pseudo-cycle. 

(ii) If $c'\in\imag\RR\setminus\CCC$ belongs to
the same connected component of $\imag\RR\setminus\CCC$ as $c$, 
$\QQQ'=\PPP_{c'}\times\III_{\omega,S^1}$, $(\rho',I')\in{\QQQ'}^{\reg}$ 
then, denoting
$\MMM'=\MMM_{\rho',I'}(\beta,c')$ and ${\EE'}^p\to\MMM'\times\Sigma^p$
the universal bundle, we can chose $E^p_T\to T$ such that there is a
diagram like (\ref{modelfinit}) and another one with
$\EE^p\to\MMM\times\Sigma^p$ replaced by 
${\EE'}^p\to\MMM'\times\Sigma^p$, and the two resulting psedo-cycles
are bordant.
\label{espseudocycle}
\end{theorem}
\begin{pf}
Recall that we write $F_1,\dots,F_r$ for the connected components 
of the fixed point set $F^f$. For any natural number
$K$, $\KK$ will denote the set $\{1,\dots,K\}\subset\NN$ and 
$\KK_0=\{0\}\cup\KK$. 
We define the framing $D=D(\Sigma^{\cusp},E,A,\phi,c)$ of the cusp 
$\rho$-THC $(\Sigma^{\cusp},E,A,\phi,c)$ to be the following set of 
data: 

\begin{enumerate}
\item The element $c\in\imag\RR$.
\item The class $\beta_0={\rho_E}_*{\phi_0}_*[\Sigma_0]\in H_2^{S^1}(F;\ZZ)$.
\item The number $K$ of bubbles in $\Sigma^{\cusp}$. 
\item Homology classes $B_1,\dots,B_K\in H_2(F;\ZZ)$ describing the
image of the bubbles $\Sigma_k$ in $\FFF$ (these classes are well
defined because the action of $S^1$ on $H_*(F;\ZZ)$ is trivial).
\item A set $S\subset\KK^2$ containing the pairs $(i,j)$
such that $\phi_i(X_i)=\phi_j(X_j)$.
\item For every $k\in\KK$ the tuple $(L_k,\Gamma_k,\rho_k)$ such that,       
after identifying $\FFF_{x_k}\simeq F$,
the bubble map $\phi_k$ belongs to 
$\MMM(L_k,\Gamma_k,\rho_k;B_k)$. 
\item For any $k\in\KK$ such that $L_k=S^1$, the number $c(k)$ such
that $\phi(\Sigma_k)\subset F_{c(k)}$;
and for any $k\in\KK$ such that $\Gamma_k\neq 1$, numbers
$c(k)_+$ and $c(k)_-$ such that $\phi(x^k_\pm)\in F_{c(k)_\pm}$, where
$x^k_{\pm}\in\Sigma_k$ are the critical points of the bubble $\Sigma_k$.
\item A set $C\subset\KK_0^2$ containing the pairs $(i,j)$
such that $i<j$ and $X_i\cap X_j\neq\emptyset$.
\begin{enumerate}
\item A partition $C=C_{00}\cup C_{01}\cup C_{10}\cup C_{11}$ defined
as follows. For any pair $(i,j)\in C$, let $x=X_i\cap X_j$.
Put $\epsilon(i)$ to be $1$ if $x\in X_i$ is a critical point with respect to
$(L_i,\Gamma_i,\rho_i)$ and $0$ otherwise (see Definition \ref{defcritic}), 
and define $\epsilon(j)$ similarly. Then $(i,j)$ belongs to 
$C_{\epsilon(i)\epsilon(j)}$.
\item A set $C''\subset\KK^3$ containing the sets $(i,j,k)$
such that $\Gamma_i(X_i\cap X_j)\cap(X_i\cap X_k)\neq\emptyset$
(this is a subset of $X_i$).
\end{enumerate}
\end{enumerate}

\begin{definition}
We will say that two cusp $\rho$-THCs 
$(\Sigma^c,E,A,\phi,c)$ and $({\Sigma'}^c,E',A',\phi',c)$ are
equivalent if
\begin{enumerate}
\item they have the same framing;
\item there is an isomorphism $g:E\stackrel{\simeq}{\longrightarrow}E'$
of bundles over $\Sigma_0=\Sigma=\Sigma_0'$ such that $g^*A'=A$ and
$g^*\phi_0'=\phi_0$;
\item let $K$ be the number of bubbles in both $\rho$-THCs; for any
$1\leq k\leq K$ there is an isomorphism $g_k:\Sigma_k\to\Sigma_k'$
such that $\phi_k=\phi_k'g_k$; furthermore, if $\Gamma_k\neq\{1\}$, 
then $g_k$ maps the critical points of $\Sigma_k$
to those of $\Sigma_k'$.
\end{enumerate}
\end{definition}

For example, if two cusp $\rho$-THCs are $\rho$-THCs (i.e., they have
no bubbles) then they are equivalent if and only if they are gauge
equivalent. 

\begin{definition}
We denote by $\MMM_{\rho,I}'(D)$ the set of equivalence classes 
of cusp $\rho$-THCs with framing $D$. 
We define the total homology class of the frame $D$
to be the equivariant homology class $\beta(D)=\beta_0+\iota_*B_1+
\dots+\iota_*B_K\in H_2^{S^1}(F;\ZZ)$. 
\label{coho}
\end{definition}

\subsubsectionr{}
The maps $\phi_k$ in a cusp curve may be multicovered. This means that
$\phi_k$ factors as $\phi_k'\circ r_k$, where $r_k:\Sigma_k\to\Sigma_k$
is a ramified covering. If $r_k$ has maximal degree, we will call
$\phi_k'$ the simplification of $\phi$. This will be a simple map.

For any cusp $\rho$-THC with frame $D$ we make the following reduction
process. First we forget the bubbles whose map to $\FFF$ is constant, 
then we substitute the bubble maps $\phi_k$ by their 
simplifications $\phi_k'$, and then we identify bubbles with the same image 
in $\FFF$. Finally, if necessary we forget some intersection points
in order that no two irreducible components of the cusp we have obtained
intersect at more than one point.
After this process we end up with another cusp $\rho$-THC with
frame $\ov{D}$. We call the resulting cusp $\rho$-THC a reduced cusp. 
We will denote $\MMM_{\rho,I}(\ov{D})\subset\MMM_{\rho,I}'(\ov{D})$
the set of equivalence classes of reduced cusp $\rho$-THCs with
framing $\ov{D}$.

Note that the total homology
class of $\ov{D}$ will not necessarily be equal to that of $D$.
If $\beta_0,\ov{B}_1,\dots,\ov{B}_K$ are the homology classes
of $\ov{D}$ we will have
$$\beta(D)=\beta_0+r_1\iota_*\ov{B}_1+\dots+r_K\iota_*\ov{B}_K,$$
where $r_k\geq 1$ are integers. This motivates the following definition.

\begin{definition}
If $\beta\in H_2(F_{S^1};\ZZ)$
and the homology classes $\beta_0,B_1,\dots,B_K$ of a frame $D$ satisfy
$\beta=\beta_0+r_1\iota_*B_1+\dots+r_K\iota_*B_K$ 
for some integers $r_k\geq 1$, then we will
say that the frame $D$ is $\beta$-admisible.
\end{definition}

We will denote $\DDD(B,c)$ the set of $\beta$-admisible framings $D$ such
that $c(D)=c$. This is obviously a numerable set. It contains a 
distinguished element $D^T$ which represents the cusp curve with no bubbling.
We will call $D^T$ the top framing of $B$.

Through all the rest of the
proof $\dim$ and $\codim$ will denote real dimension
and codimension. Recall that the dimension of $F$ is $2n$.

\subsubsectionr{}
For any tuple $(L,\Gamma,\rho)$ denote by 
$\Map_{\fibr}(L,\Gamma,\rho)$ the set of maps 
$\phi:\CP^1\to \FFF$ whose image is included in a single fibre
$\FFF_x$ and such that $\phi\in\Map(L,\Gamma,\rho)$. We have
$$\Map_{\fibr}(L,\Gamma,\rho)=E\times_{S^1}\Map(L,\Gamma,\rho).$$
Define for any homology class $B\in H_2(F;\ZZ)$
$$\MMM_{\III,\fibr}(L,\Gamma,\rho;B)
=\left\{(\phi,I)\in\Map_{\fibr}(L,\Gamma,\rho)
\times \III_{\omega,S^1}\Big|
\begin{array}{l}\ov{\partial}_I\phi=0,\ \phi_*[\Sigma]=B\\
\phi\mbox{ simple }\end{array}\right\}.$$
Incidentally, this is the moduli space used to defined
fibrewise and equivariant quantum cohomology by Givental, Kim
and Lu (see \cite{GiKm, Lu}).
Note that we have 
$\MMM_{\III,\fibr}(L,\Gamma,\rho;B)
=E\times_{S^1}\MMM_{\III}(L,\Gamma,\rho;B).$
Let $\pi:\MMM_{\III,\fibr}(L,\Gamma,\rho;B)\to\III_{\omega,S^1}$ be
the projection. Just as in the proof of Theorem \ref{modulireg} 
one can check that, for generic $I\in\III_{\omega,S^1}$, the preimage 
$\MMM_{I,\fibr}(L,\Gamma,\rho;B)=\pi^{-1}(I)$ is a smooth manifold of
dimension 
$$\dim\MMM_{I,\fibr}(L,\Gamma,\rho;B)=\dim\MMM_I(L,\Gamma,\rho;B)+2.$$

\subsubsectionr{}
Fix a $\beta$-admisible frame $D$ and
suppose that the element $c=c(D)\in\imag\RR$ lies in the complementary
of $\CCC$. 
Let $K$ be the number of bubbles,
$C\subset \{0,1,\dots,K\}^2$ the set of pairs describing which 
irreducible components intersect, $\beta_0,B_1,\dots,B_K$ the homology
classes of $D$, and $(L_k,\Gamma_k,\rho_k)$ the tuples telling the 
moduli in which $\phi_k$ sits. We denote 
$\Sigma=\Sigma_0,\Sigma_1,\dots,\Sigma_K$
the irreducible components of the cusps with frame $D$.
Define for any $k$ the group $G_k$ to be $\CC^*$ if $\Gamma_k\neq 1$ and 
$\PSL(2;\CC)$ if $\Gamma_k=1$. 
When $\Gamma_k\neq 1$, we make $G_k$ act on $\CP^1$ by
rotations keeping fixed the points $[1:0]$ and $[0:1]$, and in the
other case we consider the usual action of $G_k=\PSL(2;\CC)$ on $\CP^1$. 
The group $G_k$ acts effectively on $\MMM_I(L_k,\Gamma_k,\rho_k;B_k)$
by reparametrization: any $s\in\MMM_I(L_k,\Gamma_k,\rho_k;B_k)$ is mapped
by $g\in G_k$ to $g(s):=s\circ g:\CP^1\to F$.

\subsubsectionr{}
Let us write
$$\MMM_{\III}^*(D)=\prod_{k=1}^{K}
\MMM_{\III,\fibr}(L_k,\Gamma_k,\rho_k;B_k)\setminus\Delta,$$
where $\Delta$ is the multidiagonal, that is, the set of elements
$(s_1,\dots,s_K)$ such that $s_i=s_j$ for some $i\neq j$. 
$\MMM_{\III}^*(D)$ parameterizes
tuples of $K$ different holomorphic maps $\phi_k:\CP^1\to \FFF$
whose image is contained in any fibre.
Reasoning exactly like in the proof of Theorem 
\ref{modulireg} one proves
that $\MMM_{\III}^*(D)$ is a smooth Banach manifold.
Similarly, if $I\in\III_{\omega,S^1}$, we define
$$\MMM_I^*(D)=\prod_{k=1}^{K}
\MMM_{I,\fibr}(L_k,\Gamma_k,\rho_k;B_k)\setminus\Delta,$$
where $\Delta$ is the multidiagonal (and does not coincide with the
previous one).

\subsubsectionr{}
Let $\FFF^f=E\times_{S^1}F^f$.
For any pair $e=(j,k)\in C$ we define
$$\FFF(e)=\left\{\begin{array}{ll}
\FFF\times\FFF & \mbox{ if $e\in C_{00}$,}\\
\FFF\times\FFF^f & \mbox{ if $e\in C_{01}$,}\\
\FFF^f\times\FFF & \mbox{ if $e\in C_{10}$,}\\
\FFF^f\times\FFF^f & \mbox{ if $e\in C_{11}$,}
\end{array}\right.$$
and we write $\Delta(e)\subset\FFF(e)$ for the diagonal in $\FFF(e)$.
We also set $\Sigma(e)=\Sigma(e,j)\times\Sigma(e,k)$, where
$\Sigma(e,j)$ is defined as follows:
\begin{itemize}
\item if $j=0$ or $\Gamma_j=\{1\}$, then $\Sigma(e,j):=\CP^1$;
\item if $j>0$ and $\Gamma_j\neq \{1\}$, then
\begin{itemize}
\item if $\Sigma_j\cap\Sigma_k=x^j_+$ (resp. $x^j_-$) then
  $\Sigma(e,j)=[1:0]\in\CP^1$ (resp. $[0:1]\in\CP^1$);
\item if $\Sigma_j\cap\Sigma_k\neq x^j_\pm$ then
  $\Sigma(e,j):=\CP^1$,
\end{itemize}
\end{itemize}
and $\Sigma(e,k)$ is defined similarly.
Observe that the $\Sigma(e)$ depends on the framing $D$ (in particular,
on the isotropy pairs of the bubbles of $\Sigma^{\cusp}$), and not only
on $\Sigma^{\cusp}$.

Define
$$\bM_{\PPP,\III}(\beta,c)=
\{(A,\phi,\rho,I)\in\AAA\times\SSS\times\QQQ
\mid{\phi_E}_*[\Sigma]=\beta\text{ and $(A,\phi,\rho,I)$ 
satisfies (\ref{equs2})}\}.$$
We then have an evaluation map
$$\ev_C:\bM_{\PPP,\III}(\beta_0,c)\times\MMM_{\III}^*(D)
\times\prod_{e\in C}\Sigma(e)\to\prod_{e\in C}\FFF(e)$$
and a projection
$$\Theta:\bM_{\PPP,\III}(\beta_0,c)\times\MMM_{\III}^*(D)
\times\prod_{e\in C}\Sigma(e)\to\III_{\omega,S^1}^{K+1}.$$
Let $\Delta_{\III}$ be the diagonal in 
$\III_{\omega,S^1}^{K+1}$. Since $c\in\imag\RR\setminus\CCC$, Theorems
\ref{modulipert} and \ref{modulireg} imply that 
$\Theta^{-1}(\Delta_{\III})$ is a smooth
Banach manifold (this would not be true if we had not removed the
multidiagonal in the definition of $\MMM_{\III}^*(D)$ --- a similar
thing occurs in Lemma 4.9 in \cite{RuTi}).
Let us define
$$\bR_{\PPP,\III}(D)=\Theta^{-1}(\Delta_{\III})
\cap\ev_C^{-1}(\prod\Delta(e)).$$
Lemma \ref{submersioeq}, together with its 
nonequivariant version Lemma 6.1.2 in
\cite{McDS1}, implies that the restriction of $\ev_C$ to 
$\Theta^{-1}(\Delta_{\III})$ is a submersion. Hence, 
$\bR_{\PPP,\III}(D)$ is a smooth Banach manifold.
On the other hand, the gauge group $\GGG$ acts freely on
$\bR_{\PPP,\III}(D)$ (because it acts freely on $\bM(\beta_0,c)$), so
$$\RRR_{\PPP,\III}(D):=\bR_{\PPP,\III}(D)/\GGG$$
is also a Banach manifold.

Consider the projection 
$$q:\RRR_{\PPP,\III}(D)\to\QQQ=\PPP_c\times\III_{\omega,S^1}.$$
This is a Fredholm map, to which we may apply Sard--Smale theorem. We
deduce that there exists a subset $\QQQ^{\reg}(D)\subset\QQQ$ of Baire of
the second category such that for any $(\rho,I)\in\QQQ^{\reg}(D)$ the
set $$\RRR_{\rho,I}:=q^{-1}(\rho,I)$$ is a smooth Banach manifold of
dimension equal to the index of the differential of $q$, i.e., to the
virtual dimension.

\begin{lemma}
Suppose that $\DDD(\beta,c)\ni 
D\neq D^T$. If $(\rho,I)\in\QQQ^{\reg}(D)$ then
$$\dim \RRR_{\rho,I}(D)\leq
\dim \MMM(\beta,c)-2+\sum_{k=1}^K\dim G_k.$$
\label{boundcodim}
\end{lemma}
\begin{pf}
Let $\bR_{\rho,I}(D)$ be the preimage of $\RRR_{\rho,I}(D)$ by the
projection map $\bR_{\PPP,\III}(D)\to\RRR_{\PPP,\III}(D)$. We then
have a commutative diagram
$$\xymatrix{\bR_{\rho,I}(D) \ar[r]^-{\tilde{\iota}}\ar[d] &
\bM_{\rho,I}(\beta_0,c)\times\MMM^*(D)\times\prod\Sigma(e)\ar[d] \\
\RRR_{\rho,I}(D)\ar[r]^-{\iota} &
(\bM_{\rho,I}(\beta_0,c)\times\MMM^*(D)\times\prod\Sigma(e))/\GGG.}$$
Since $\GGG$ acts smoothly and freely on the top rows and the vertical
arrows are projection to $\GGG$ orbits, it follows that the
codimension of $\iota$ is the same as the codimension of
$\tilde{\iota}$. Now, $\bR_{\rho,I}(D)$ is the preimage of
$\prod\Delta(e)$ by the evaluation map
$$\ev_C(\rho,I):\bM_{\rho,I}(\beta_0,c)\times
\MMM_I^*(D)\times\prod\Sigma(e)\to\prod\FFF(e),$$  
and since $(\rho,I)\in\QQQ^{\reg}(D)$, it follows that this is a
submersion. Consequently, the codimension of $\tilde{\iota}$ is equal
to the codimension of the inclusion
$\prod\Delta(e)\subset\prod\FFF(e)$. Hence,
\begin{align*}
\dim\RRR_{\rho,I}(D)-\sum\dim G_k &=
\dim\MMM_{\rho,I}(\beta_0,c)+\dim\MMM^*_I(D)+\sum\dim\Sigma(e) \\
&+ \sum\dim\Delta(e)-\sum\dim\FFF(e)-\sum\dim G_k.
\end{align*}

To bound this dimension we divide the set of bubbles $\KK$ in 
three subsets. Let $\SS$ (resp. $\TT$ and $\UU$) denote the set
of $k\in\KK$ such that $L_k=1$, $\Gamma_k=1$ (resp.
$L_k=1$, $\Gamma_k\neq 1$ and $L_k=S^1$, $\Gamma_k=1$).
Let $S=|\SS|$, $T=|\TT|$ and $U=|\UU|$. Theorems \ref{indtor} and
\ref{inds} imply the following.
\begin{itemize}
\item If $k\in\SS$ then $\dim G_k=6$ and
\begin{align*}
\dim\MMM_{I,\fibr}(L_k,\Gamma_k,\rho_k;B_k)
&= 2+2\la c_1(TF),B_k\ra+2n \\
&= 2+2\la c_1^{S^1}(TF),\iota_*B_k\ra+2n.
\end{align*}
\item If $k\in\TT$ then $\dim G_k=2$ and
\begin{align*}
\dim\MMM_{I,\fibr}(L_k,\Gamma_k,\rho_k;B_k)
&\leq 2+2\la c_1(TF),B_k\ra+2n-4 \\
&= 2\la c_1^{S^1}(TF),\iota_*B_k\ra+2n-2,
\end{align*}
by Lemma \ref{destor} and Condition (\ref{cond3}).
\item If $k\in\UU$ then $\dim G_k=6$ and
\begin{align*}
\dim\MMM_{I,\fibr}(L_k,\Gamma_k,\rho_k;B_k) 
&=2+2\la c_1(TF^f),B_k\ra+\dim F_{c(k)} \\
&\leq2+2\la c_1(TF),B_k\ra+\dim F_{c(k)} \\
&=2+2\la c_1^{S^1}(TF),\iota_*B_k\ra+\dim F_{c(k)},
\end{align*}
by Condition (\ref{cond2}).
\end{itemize}
On the other hand, since $D$ is $B$-admissible, using Condition
(\ref{cond1}) we obtain 
$$\la c_1^{S^1}(TF),\beta_0\ra+\sum_{k=1}^K\la c_1^{S^1}(TF),\iota_*B_k\ra
\leq\la c_1^{S^1}(TF),\beta\ra.$$
Hence, 
\begin{align*}
\dim\RRR_{\rho,I}(D)-\sum\dim G_k 
&\leq 2\la c_1^{S^1}(TF),\beta\ra+2(n-1)(1-g) \\
&+(S+T)2n+\sum_{k\in\UU}\dim F_{c(k)}-4K \\
&+\sum_{e\in C}
\dim\Sigma(e)+\dim \Delta(e)-\dim\FFF(e).
\end{align*}
To find an upper bound for the last two terms we proceed as follows.
Since for any $e$ we have 
$$\dim\Sigma(e)+\dim \Delta(e)-\dim\FFF(e)\leq 0,$$ an
upper bound for $\sum_{e\in C'}$ where $C'\subset C$ will also give a
bound on $\dim\RRR_{I,\sigma}(D)$. So we take any subset $C'\subset C$ of $K$
elements with the following property. The graph whose vertices are the 
elements of $\KK_0$ and which has an edge joining $i$ to $j$ if either 
$(i,j)$ or $(j,i)$ belong to $C'$ is connected. This implies that 
$C'\nsubseteq C_{11}$ (because otherwise the vertex $0\in\KK_0$ would be 
disconnected from the rest). Take an injective map
$$v:C'\to\KK$$
which assigns to $(i,j)$ either $i$ or $j$. Let $k\in\KK$ and $e=v^{-1}(k)$.
\begin{itemize}
\item If $k\in\SS$ then $\dim\Sigma(e)+\dim\Delta(e)-\dim\FFF(e)=-2n+2$.
\item If $k\in\TT$ and $e\notin C_{11}$ then
$\dim\Sigma(e)+\dim\Delta(e)-\dim\FFF(e)\leq-2n+2$ and if $e\in C_{11}$ then 
$$\dim\Sigma(e)+\dim\Delta(e)-\dim\FFF(e)\leq
\max\{-\dim F_{c(k)_+},-\dim F_{c(k)_+}\}-2\leq -2n+4,$$
by Condition (\ref{cond3}).
\item If $k\in\UU$ then 
$\dim\Sigma(e)+\dim\Delta(e)-\dim\FFF(e)=-\dim F_{c(k)}+2$.
\end{itemize}
Since $C'\nsubseteq C_{11}$, we get
\begin{align*}
\dim\RRR_{\rho,I}(D)-\sum\dim G_k
&\leq 2\la c_1^{S^1}(TF),B\ra+2(n-1)(1-g)-2 \\
&=\dim\MMM_{\rho,I}(B,c)-2,
\end{align*}
which is what we wanted to prove.
\end{pf}

On the other hand, note that $\bR_{\rho,I}(D)$ is invariant under the
action of $\prod G_k$ in
$\bM_{\rho,I}(\beta_0,c)\times\MMM_I^*(D)\times\prod\Sigma(e)$. This
action is free and commutes with the action of $\GGG$ (indeed, $\prod
G_k$ acts on $\MMM^*_I(D)$ by reparametrization of bubbles, while
$\GGG$ acts on it by acting on its image ---hence the first action is on
the right and the second one is on the left). So $\RRR_{\rho,I}(D)$
inherits an action of $\prod G_k$, and it is not difficult to see that
this action is free. 

Let us define $$\QQQ^{\reg}=\bigcap_{D\in\DDD(\beta,c)}
\QQQ^{\reg}(D),$$ and assume that $(\rho,I)\in\QQQ^{\reg}$.
Just as in Subsection \ref{findimtarget} we have, for any framing
$D\neq D^T$, an evaluation map taking values in the compact smooth model of
the universal bundle
$$\ev^0_D:\RRR_{\rho,I}(D)\times(\bigcup\Sigma_k)^p\to
E_T^p\times_{(S^1)^p}F^p.$$
If we make $\prod G_k$ act on $(\bigcup\Sigma_k)^p$ as well as on
$\RRR_{\rho,I}$ then the map $\ev^0_D$ is invariant, so it descends
to a map 
$$\ev_D:(\RRR_{\rho,I}(D)\times(\bigcup\Sigma_k)^p)/\prod G_k\to
E_T^p\times_{(S^1)^p}F^p.$$
By Lemma \ref{boundcodim}, the domain of this map has dimension
$$\dim(\RRR_{\rho,I}(D)\times(\bigcup\Sigma_k)^p)/\prod G_k
\leq\dim(\MMM(\beta,c)\times(\bigcup\Sigma_k)^p)-2.$$
Finally, Theorem \ref{compactificacio} on compactness tells us that
$$\overline{s^T_{F^p}\Phi(\MMM(\beta,c)\times\Sigma^p)}
\subset s^T_{F^p}\Phi(\MMM(\beta,c)\times\Sigma^p)
\cup\bigcup_{D^T\neq D\in\DDD(\beta,c)}\Im \ev_D.$$
Indeed, $\RRR_{\rho,I}(D)/\prod G_k$ parametrizes the set
$\MMM_{\rho,I}(D)$ of equivalence
classes of cusp $\rho$-THCs with framing $D$.
This finishes the proof that $s^T_{F^p}\Phi$ is a
pseudo-cycle. 

The second part of the theorem is proved following exactly the same
method. 
\end{pf}

\subsection{Definition of the invariants}
Proceeding as in Subsection \ref{toprest} with $P=E^p_T$, $K=(S^1)^p$
and $V=F^p$, we get a map
$$c_T^p:E_T^p\times_{(S^1)^p}F^p\to (F^p)_{(S^1)^p}.$$
On the other hand, there is a projection map
$$\nu_T:T=(\AAA_{\Fl}\times S)/\GGG\times\Sigma^p\to
\AAA_{\Fl}/\GGG=Jac_d(\Sigma).$$
Note that $(ES^1)^p$ is a contractible space on which $(S^1)^p$ acts
freely. Hence
\begin{equation}
(F^p)_{(S^1)^p}=((ES^1)^p\times F^p)/(S^1)^p=
((ES^1\times F)/S^1)^p=(F_{S^1})^p.
\label{isocoheq}
\end{equation}
Finally, let
$$\psi_{\Sigma,\beta,c}:H_*(E^p_T\times_{(S^1)^p}F^p;\ZZ)\to\ZZ$$
be the map induced by the pseudo-cycle $s^T_{F^p}$ using Lemma
\ref{aplihom}. 

Now let $\alpha_1,\dots,\alpha_p\in H_{S^1}^*(F;\ZZ)$. Combining
K{\"u}nneth with the isomorphism (\ref{isocoheq}) we may view
$$\alpha_1\otimes\dots\otimes\alpha_p\in H^*_{(S^1)^p}(F^p;\ZZ).$$
Let also $$\gamma\in H^*(\AAA/\GGG;\ZZ)=H^*(\Jac_d(\Sigma);\ZZ).$$
We define the Hamiltonian Gromov--Witten invariant of
$\alpha_1,\dots,\alpha_p,\gamma$ at $\beta,c$ to be
$$\Psi_{\Sigma,\beta,c}(\alpha_1,\dots,\alpha_p,\gamma):=
\psi_{\Sigma,\beta,c}(PD({c_T^p}^*(\alpha_1\otimes\dots\otimes\alpha_p)
\cup\nu_T^*\gamma)),$$ 
where $PD$ denotes Poincar{\'e} Dual.
It follows from Lemma \ref{aplihom} and Theorem \ref{espseudocycle} that 
$\Psi_{\Sigma,\beta,c}(\alpha_1,\dots,\alpha_p,\gamma)$ only depends on the
deformation class of $\omega$ as a $S^1$-invariant symplectic
structure, on he connected component in $\imag\RR\setminus\CCC$ in
which $c$ lies, the (co)homology classes
$\beta,\alpha_1,\dots,\alpha_p,\gamma$, and the Riemann surface
$\Sigma$ (in fact, it only depends on the genus of $\Sigma$ and not on
its conformal class, as can be
proved using the same cobordism methods as above). By Subsection
\ref{sb:regularity} 
the invariants are independent of the Sobolev norm which was used 
to complete $\AAA$, $\SSS$ and $\GGG$.

\end{document}